\documentclass[11pt,a4paper]{article}
\usepackage[T1]{fontenc}
\usepackage[utf8]{inputenc}
\usepackage[british]{babel}
\usepackage{lmodern,microtype}
\usepackage[margin=15mm]{geometry}
\usepackage{amsmath,amssymb,amsthm,mathtools}
\usepackage{booktabs,array,mathrsfs}
\usepackage{xcolor}
\usepackage{tikz}
\usepackage{tikz-cd}
\usetikzlibrary{arrows.meta,calc,positioning,decorations.pathreplacing,decorations.markings}
\usepackage{cite}
\usepackage[hidelinks]{hyperref}
\usepackage{enumitem}
\usepackage[font=small,labelfont=bf]{caption}
\hypersetup{pdftitle={Homological reduction and entropy of discontinuous interval maps},pdfauthor={Henrique M. Oliveira}}
\newtheorem{theorem}{Theorem}[section]
\newtheorem{proposition}[theorem]{Proposition}
\newtheorem{corollary}[theorem]{Corollary}
\newtheorem{lemma}[theorem]{Lemma}
\theoremstyle{definition}

\theoremstyle{remark}
\newtheorem{remark}[theorem]{Remark}
\DeclareMathOperator{\imageop}{im}
\DeclareMathOperator{\spanop}{span}

\newcommand{\Id}{\mathrm{Id}}
\newcommand{\eps}{\varepsilon}
\newcommand{\calB}{\mathcal{B}}
\newcommand{\calG}{\mathcal{G}}
\newcommand{\calP}{\mathcal{P}}
\newcommand{\R}{\mathbb R}

\newcommand{\ehat}{\widehat F}
\newcommand{\sig}{\operatorname{sgn}}
\tikzset{
 map branch/.style={black,line width=1.1pt},
 map frame/.style={black,line width=.7pt},
 map axis/.style={black!30,line width=.35pt},
 open end/.style={circle,draw=black,fill=white,line width=.65pt,
                  inner sep=0pt,minimum size=3.8pt},
 closed end/.style={circle,fill=black,inner sep=0pt,minimum size=3.8pt},
 jump edge/.style={line width=.8pt,postaction={decorate},
  decoration={markings,mark=at position .78 with
   {\arrow{Stealth[length=2mm,width=1.4mm]}}}},
 edge label/.style={fill=white,inner sep=1.5pt,font=\small}
}
\newcommand{\one}{\mathbf1}
\title{Homological reduction and entropy\\of discontinuous interval maps}
\author{Henrique M. Oliveira\\[3pt]
\small Centro de An\'alise Matem\'atica, Geometria e Sistemas Din\^amicos (CAMGSD),\\
\small Departamento de Matem\'atica, Instituto Superior T\'ecnico, Universidade de Lisboa,\\
\small Av. Rovisco Pais 1, 1049--001 Lisbon, Portugal\\[2pt]
\small \texttt{henrique.m.oliveira@tecnico.ulisboa.pt}\\[2pt]
\small ORCID: 0000-0002-3346-4915}
\date{}
\begin{document}
\maketitle
\begin{abstract}
We develop a homological approach to the topological entropy of discontinuous interval maps with arbitrary branch orientations. Relations among the one-sided images of discontinuities reduce the kneading determinant, while symmetries yield further simplifications. We determine the smallest possible dimension of this reduction for each orientation pattern and show that this minimum can coexist with maximal entropy. The general theory is applied to maps with increasing or decreasing branches, extending earlier results for bimodal maps and yielding exact formulas for entropy and periodic-orbit counts.
\end{abstract}
\section{Introduction}\label{sec:intro}
We study how relations among the one-sided images of discontinuities reduce the kneading determinant of an interval map. Our main result determines the minimum dimension of this reduction for every prescribed branch-orientation pattern. This minimum never exceeds three and can coexist with maximal entropy for each fixed number of branches. We also obtain exact cylinder-counting formulas and identify the determinant factor carrying the entropy under compatible orbital symmetries.
The motivating examples come from the complex tangent family. For a purely imaginary parameter $i\beta$, its second iterate on the real axis is
\[
 f_\beta(x)=-\beta\tanh(\beta\tan x).
\]
For $\pi/2<\beta<3\pi/2$, its restriction to $[-\beta,\beta]$ has two discontinuities and three decreasing branches. In joint work with Alves and Sousa Ramos \cite{AlvesOliveiraRamos2004}, a kneading determinant was used to calculate its entropy and factor its dynamical zeta function. That construction already grouped coincident one-sided image values. In \cite{Oliveira2015}, I developed the symbolic and homological description of odd discontinuous bimodal maps and related the kneading determinant to the spectrum of a finite Markov matrix.

These examples combine two features: the jumps share image data, and reflection relates their subsequent orbits. The first feature suggests reducing the determinant before computing the itineraries; the second suggests a further factorisation. With more discontinuities, shared images need not come from symmetry. With mixed branch orientations, the signs attached to the two sides of a cut also change. The present construction separates these effects by recording an image point together with its side of approach, and then taking the image of the signed jump boundary.

Entropy for interval maps is described by orbit separation and itinerary growth \cite{AdlerKonheimMcAndrew1965,Bowen1971,Dinaburg1971,MisiurewiczSzlenk1980,LindMarcus1995,AlsedLlibreMisiurewicz2000}. Milnor and Thurston \cite{MilnorThurston1988} encode critical-orbit data in a kneading determinant; Hofbauer \cite{Hofbauer1981} analyses the symbolic structure of piecewise continuous maps. Preston \cite{Preston1989}, Baladi and Keller \cite{BaladiKeller1990}, Baladi and Ruelle \cite{BaladiRuelle1994}, and Baladi \cite{Baladi1995} develop determinant and transfer-operator methods in the discontinuous setting. Rugh and Tan \cite{RughTan2015} establish weighted kneading identities for arbitrary orientations. We use their unweighted entropy theorem; the reduction concerns the dimension and structure of its interior determinant.

Homological approaches include the bimodal construction of Lampreia and Sousa Ramos \cite{LampreiaRamos1997}, the functorial theory of Alves and Sousa Ramos \cite{AlvesRamos1999,AlvesRamos2004}, and the graph-map results of Alves, Hric and Sousa Ramos \cite{AlvesHricRamos2005}. Silva and Sousa Ramos \cite{SilvaRamos2002,SilvaRamos2003} treat Lorenz renormalisation and all four two-branch orientation types. Further discontinuous kneading constructions appear in \cite{RochaRamos2003,RochaRamos2004}, while entropy conventions, matching, and systems with holes are studied in \cite{CalderonVillar2024,BruinCarminatiMarmiProfeti2019,SunLiDing2024,DingSteinerSun2026}.

The signed boundary identifies the directions on which the kneading operator is the identity (Theorem~\ref{thm:reduction}). Its rank is computed by the bipartite components of the jump graph (Theorem~\ref{thm:rank}), using classical incidence theory \cite{Zaslavsky1982,CvetkovicRowlinsonSimic2007}. The minimum rank for a prescribed orientation word equals the number of side types represented by the image-germ pairs: both approached from the right, both from the left, or one from each side. In particular, the adjacent branch-orientation pairs $(+,+)$ and $(-,-)$ belong to the same third class. For $n$ branches, an affine Markov construction attains this minimum and has entropy within $\log(n/(n-1))$ of $\log n$ (Theorems~\ref{thm:sharp-pattern-rank} and~\ref{thm:sharp-high-entropy}). A further construction preserves the minimal rank and admits every finite branch word (Proposition~\ref{prop:full-entropy-minimal-rank}).

The cylinder identity in Proposition~\ref{prop:counting} supplies the counting formulas. For maps with positive entropy, it also identifies the symmetry factor whose first positive zero determines the entropy (Corollary~\ref{cor:symmetric-entropy}). For a finite Markov partition, Theorem~\ref{thm:Markov} retains the endpoint and germ-cycle factors; related comparisons occur in \cite{GlendinningHall1996,Toulson2019}. Sections~\ref{sec:homology}--\ref{sec:sharp} develop the general reduction, entropy formulas, symmetry and sharp realisations. Section~\ref{sec:mixed} illustrates mixed orientations.

The final applications treat constant orientations. Section~\ref{sec:decreasing} recovers the odd bimodal formula of \cite{Oliveira2015} and gives decreasing families with exact entropy and cylinder counts. Section~\ref{sec:increasing} treats increasing branches. At a finite critical return of a symmetric Lorenz map, the antisymmetric Markov factor is cyclotomic and the symmetric factor is determined by the critical itinerary. For uniformly expanding odd Lorenz maps, sign reversal preserves the number of cycles of every least period at least three, even without a finite critical orbit. These statements use the increasing--decreasing correspondence studied in \cite{AnusicBruinCinc2020}; the connection of constant-slope factors with Lorenz flows is developed in \cite{CholewaIgra2026}. Section~\ref{sec:shadowing} constructs a surjective discontinuous interval map with shadowing, answering the existence question in \cite{BuraGoodSamuel2026}.

\section{Maps, germs and signed itineraries}\label{sec:setting}

\subsection{Branches and entropy}

Fix $m\geq1$, $I=[a,b]$, and
\[
 a=c_0<c_1<\cdots<c_m<c_{m+1}=b,\qquad
 J_j=(c_j,c_{j+1})\quad(0\leq j\leq m).
\]
For each $0\leq j\leq m$, let $F_j:\overline{J_j}\to I$ be continuous and strictly monotone. Its restriction to $J_j$ is the $j$th branch. Let
\[
 \eps_j=\begin{cases}
 +1,&F_j\text{ is increasing},\\
 -1,&F_j\text{ is decreasing}.
 \end{cases}
\]
The tuple $\eps=(\eps_0,\ldots,\eps_m)$ is the \emph{orientation word}. A cut $c_i$ is a \emph{genuine jump} when $F_{i-1}(c_i)\ne F_i(c_i)$. The cuts are genuine jumps unless a continuous monotonicity cut is explicitly admitted; the recurrence and determinant identities below apply in either case. Write $F$ for the map given by the branches on $\bigcup_{j=0}^m J_j$, and $F^k$ for its $k$th iterate wherever this is defined, with $F^0$ the identity.

For $n\geq1$ and a word $w=(i_0,\ldots,i_{n-1})\in\{0,\ldots,m\}^n$, its \emph{cylinder} is
\[
 C(w)=\{x\in(a,b):F^k(x)\in J_{i_k}\text{ for }0\leq k<n\}.
\]
The empty word has cylinder $(a,b)$. A word is \emph{admissible} when its cylinder is nonempty. Each such cylinder is an open interval on which $F^k$, $0\leq k<n$, is defined and strictly monotone. Let $\ell_n(F)$ be the number of admissible words of length $n$, so $\ell_0(F)=1$. The \emph{topological entropy} used throughout is
\begin{equation}\label{eq:entropy-definition}
 h_{\mathrm{top}}(F)=\lim_{n\to\infty}\frac{1}{n}\log\ell_n(F).
\end{equation}
Write $\gamma_F=\exp(h_{\mathrm{top}}(F))$ for the cylinder growth number. The limit exists because $\ell_{n+r}(F)\leq\ell_n(F)\ell_r(F)$. In particular, $0\leq h_{\mathrm{top}}(F)\leq\log(m+1)$.

For a continuous map $T$ of a compact space, the Adler--Konheim--McAndrew entropy is
\[
 h(T)=\sup_{\mathcal U}\lim_{n\to\infty}\frac1n
 \log\nu\!\left(\bigvee_{k=0}^{n-1}T^{-k}\mathcal U\right),
\]
where $\mathcal U$ ranges over finite open covers, $\nu$ is the least cardinality of a subcover, and $\bigvee$ denotes their common refinement \cite{AdlerKonheimMcAndrew1965}. Let $\Omega_F$ be the set of points whose entire forward orbit lies in $\bigcup_j J_j$, and put $\mathbb N_0=\{0,1,\ldots\}$. The \emph{itinerary} of $x\in\Omega_F$ is the sequence $(i_k)_{k\in\mathbb N_0}$ defined by $F^k(x)\in J_{i_k}$. Give the alphabet $\{0,\ldots,m\}$ the discrete topology, and let $\Sigma_F$ be the closure of these itineraries in the product topology on $\{0,\ldots,m\}^{\mathbb N_0}$. The complement of $\Omega_F$ in $I$ is countable, so every nonempty cylinder contains a point of $\Omega_F$. The left shift $\sigma$ is continuous on $\Sigma_F$. Its cylinder partitions refine every finite open cover, giving $h(\sigma|_{\Sigma_F})=h_{\mathrm{top}}(F)$ by \eqref{eq:entropy-definition} \cite[Chapter~4]{LindMarcus1995}.

For $x,y\in\Omega_F$, put
\[
 d_n(x,y)=\max_{0\leq k<n}|F^k(x)-F^k(y)|.
\]
A finite subset of $\Omega_F$ is \emph{$(n,\eta)$-separated} if distinct points have $d_n$-distance greater than $\eta$; denote the supremum of its possible cardinalities by $s_n(\eta)$.

\begin{proposition}[Metric entropy]\label{prop:entropy-models}
For the maps of this section,
\begin{equation}\label{eq:entropy-bowen}
 h_{\mathrm{top}}(F)
 =\lim_{\eta\downarrow0}\limsup_{n\to\infty}
       \frac{1}{n}\log s_n(\eta).
\end{equation}
\end{proposition}
\begin{proof}
For genuine jumps, the extended branch graphs are pairwise disjoint compact subsets of $I^2$. Let $\delta>0$ be their minimum distance in the maximum metric. Points in different $n$-cylinders first differ in branch address at some time $k<n$; either their positions at time $k$ or their images at time $k+1$ differ by at least $\delta$. Choosing one point of $\Omega_F$ per cylinder gives
\[
 s_{n+1}(\delta/2)\geq\ell_n(F).
\]
Within an $n$-cylinder, all $F^k$, $0\leq k<n$, are monotone. For an ordered separated set, each adjacent pair contributes more than $\eta$ to at least one coordinate variation. Their total variation is at most $n(b-a)$, so
\[
 s_n(\eta)\leq\ell_n(F)\bigl(1+n(b-a)/\eta\bigr).
\]
Taking growth rates proves \eqref{eq:entropy-bowen}. If continuous monotonicity cuts are admitted, the same equality follows from \cite[Theorems~2.4 and 2.6]{CalderonVillar2024}.
\end{proof}

For genuine jumps, extend $F$ by its inward endpoint limits and arbitrary cut values. The set of points whose forward orbits avoid the cuts then differs from $\Omega_F$ by at most the two endpoints, so the separated-set growth rate is unchanged. This also equals the piecewise-continuous open-cover entropy of \cite[Theorem~2.4]{CalderonVillar2024}. Values at the cuts do not enter \eqref{eq:entropy-definition}.

\subsection{One-sided orbits}

A \emph{germ} $x^\sigma$, $\sigma\in\{-1,+1\}$, records approach to $x$ from the indicated side; the endpoints have only $a^+$ and $b^-$. Germs with different base points inherit the order of those points, and at an interior point we use $x^-<x<x^+$. A germ belongs to branch $j$ when points approaching its base from its indicated side lie in $J_j$. For such a germ $g=x^\sigma$, write $F(x^\sigma)=F_j(x)$ for the one-sided image value and define
\[
 \iota(g)=j,\qquad \eps(g)=\eps_j,\qquad
 \ehat(g)=F(x^\sigma)^{\eps_j\sigma}.
\]
Thus $\ehat$ follows one-sided images, including when an image lies at a cut. The \emph{itinerary} of a germ $g$ is the sequence $(\iota(\ehat^k g))_{k\in\mathbb N_0}$. Write $\delta(x^\sigma)=\sigma$ for the side of a germ.

For $1\leq i\leq m$, put
\[
 g_i^-=\ehat(c_i^-),\qquad g_i^+=\ehat(c_i^+),
 \qquad \calG=\{g_i^-,g_i^+:1\leq i\leq m\}.
\]
Repeated germs are included only once in $\calG$. Let $E=\R^m$, with basis $e_1,\ldots,e_m$, and $V=\R^{\calG}$, with basis $[g]$ indexed by the distinct germs.

\subsection{Signed itinerary series}

For each cut $c_i$ define the centred indicator
\[
 \varphi_i(g)=\begin{cases}
 -\tfrac12,&g<c_i,\\
 +\tfrac12,&g>c_i.
 \end{cases}
 \qquad \varphi(g)=(\varphi_1(g),\ldots,\varphi_m(g))^{\mathsf T}.
\]
The order on germs makes this unambiguous even at $c_i^\pm$. The cumulative orientation during the first $n$ steps is
\[
 \eps^{(0)}(g)=1,\qquad
 \eps^{(n)}(g)=\prod_{k=0}^{n-1}\eps(\ehat^k g),
 \qquad n\geq1.
\]
For a real vector space $L$, write $L[[t]]$ for formal power series in $t$ with coefficients in $L$; all identities are initially coefficientwise. The \emph{signed itinerary series} is
\begin{equation}\label{eq:theta}
 \theta(g,t)=\sum_{n\geq0}t^n\eps^{(n)}(g)\varphi(\ehat^n g)
 \in E[[t]].
\end{equation}
In particular,
\begin{equation}\label{eq:recurrence}
 \theta(g,t)=\varphi(g)+t\eps(g)\theta(\ehat g,t).
\end{equation}
The linear map $K_t:V\to E[[t]]$ is specified by $K_t[g]=\theta(g,t)$.

\section{Signed jump homology and determinant reduction}\label{sec:homology}

\subsection{The boundary complex}

Define the signed jump boundary $\partial_\eps:E\to V$ by
\begin{equation}\label{eq:boundary}
 \partial_\eps e_i=\eps_i[g_i^+]-\eps_{i-1}[g_i^-].
\end{equation}
Regard $E$ and $V$ as chain spaces in degrees one and zero. The two-term chain complex $0\to E\xrightarrow{\partial_\eps}V\to0$ has homology
\[
 H_1^\eps=\ker\partial_\eps,\qquad
 H_0^\eps=V/\imageop\partial_\eps.
\]
Thus a degree-one cycle is a linear combination of cuts whose signed image contributions cancel. Write
\[
 \calB_\eps=\imageop\partial_\eps,\qquad
 r_\eps=\dim\calB_\eps.
\]

The interior kneading operator is $N(t):E[[t]]\to E[[t]]$ defined on basis vectors by
\[
 N(t)e_i=\theta(c_i^+,t)-\theta(c_i^-,t),
 \qquad D_F(t)=\det N(t).
\]
All operators are extended $\R[[t]]$-linearly to formal series; restrictions to subspaces use the same convention. Determinants are taken in $\R[[t]]$.

\begin{lemma}[Signed kneading recurrence]\label{lem:recurrence}
For every orientation pattern,
\begin{equation}\label{eq:factor}
 N(t)=\Id_E+tK_t\partial_\eps.
\end{equation}
\end{lemma}
\begin{proof}
Apply \eqref{eq:recurrence} to $c_i^+$ and $c_i^-$. Their centred-indicator vectors differ by $e_i$, and their first-step orientation factors are $\eps_i$ and $\eps_{i-1}$. Consequently
\[
 \theta(c_i^+,t)-\theta(c_i^-,t)
 =e_i+t\bigl(\eps_i\theta(g_i^+,t)-\eps_{i-1}\theta(g_i^-,t)\bigr),
\]
which is \eqref{eq:factor} in column $i$.
\end{proof}

\subsection{Reduction to the boundary space}

Since $\partial_\eps K_t$ takes values in $\calB_\eps[[t]]$, the operator $\Id_V+t\partial_\eps K_t$ preserves this space. We may therefore define the reduced operator
\[
 R_\eps(t)=\left(\Id_V+t\partial_\eps K_t\right)|_{\calB_\eps}.
\]
Both $N(0)$ and $R_\eps(0)$ are identity operators, so $N(t)$ and $R_\eps(t)$ are invertible over $\R[[t]]$.

The boundary map gives the exact sequence
\begin{equation}\label{eq:boundary-exact}
 0\longrightarrow H_1^\eps\longrightarrow E
 \xrightarrow{\partial_\eps}\calB_\eps\longrightarrow0
\end{equation}
where exactness means that the image of each arrow is the kernel of the next. A sequence is \emph{equivariant} for operators on its spaces when those operators commute with its arrows. We use this terminology also after extension to formal series.

\begin{theorem}[General homological reduction]\label{thm:reduction}
Sequence~\eqref{eq:boundary-exact} is equivariant for $\Id_{H_1^\eps}$, $N(t)$ and $R_\eps(t)$. In particular,
\begin{equation}\label{eq:red}
  D_F(t)=\det R_\eps(t).
\end{equation}
Thus the determinant has a representation of order $r_\eps\leq m$, independently of the number of distinct forward itineraries of the cut germs.
\end{theorem}
\begin{proof}
Equation \eqref{eq:factor} implies that $N(t)$ acts as the identity on $\ker\partial_\eps$. Furthermore,
\[
 \partial_\eps N(t)=
 \bigl(\Id_V+t\partial_\eps K_t\bigr)\partial_\eps
 =R_\eps(t)\partial_\eps.
\]
Choose a basis of $E$ beginning with a basis of $H_1^\eps$ and project the remaining basis vectors to $\calB_\eps$. The matrix of $N(t)$ is block upper triangular, with diagonal blocks $\Id_{H_1^\eps}$ and $R_\eps(t)$. Taking determinants gives \eqref{eq:red}.
\end{proof}

A \emph{chain map} on the two-term complex is a pair $(A,B)$ with $\partial_\eps A=B\partial_\eps$. It is \emph{chain-homotopic to the identity} if a map $H:V\to E$ satisfies $A-\Id_E=H\partial_\eps$ and $B-\Id_V=\partial_\eps H$. Here $\widetilde R_\eps(t)=\Id_V+t\partial_\eps K_t$ and $H_t=tK_t$ give such a homotopy \cite[Section~2.1]{Hatcher2002}. Then
\[
 N(t)-\Id_E=H_t\partial_\eps,\qquad
 \widetilde R_\eps(t)-\Id_V=\partial_\eps H_t.
\]
The induced actions on $H_1^\eps$ and $H_0^\eps$ are identities.

For computation, choose an isomorphism $U:\R^{r_\eps}\to\calB_\eps\subset V$ and put $L=U^{-1}\partial_\eps:E\to\R^{r_\eps}$. Then $\partial_\eps=UL$ and
\begin{equation}\label{eq:small-matrix}
 D_F(t)=\det\bigl(\Id_{r_\eps}+tLK_tU\bigr).
\end{equation}
A change of boundary basis conjugates this matrix. This is the finite-rank determinant reduction of \cite[Proposition~1.2]{AlvesRamos1999}, with the reduced space specified by the jump images.

\subsection{The jump graph and its rank}

The jump graph $\Gamma_F$ has vertex set $\calG$ and one undirected edge joining $g_i^-$ to $g_i^+$ for each cut $c_i$; parallel edges are retained. A loop is possible if a continuous turning point is included as a cut, but cannot occur at a genuine jump. Let $q=|\calG|$, and let $b(\Gamma_F)$ denote the number of connected components of $\Gamma_F$ that are bipartite. A component is bipartite when its vertices can be partitioned into two sets with every edge joining different sets. Write $c(\Gamma_F)$ for the total number of connected components. The following formula applies the incidence-rank theorem for signed graphs \cite{Zaslavsky1982,CvetkovicRowlinsonSimic2007} to the jump boundary.

\begin{theorem}[Jump-rank formula]\label{thm:rank}
For arbitrary branch orientations,
\begin{equation}\label{eq:rank}
  r_\eps=q-b(\Gamma_F),\qquad
 \dim H_1^\eps=m-q+b(\Gamma_F).
\end{equation}
If every branch is increasing, or every branch is decreasing, all components are bipartite and $r_\eps=q-c(\Gamma_F)$.
\end{theorem}
\begin{proof}
By the definition of germ transport,
\[
 \delta(g_i^+)=\eps_i,\qquad
 \delta(g_i^-)=-\eps_{i-1}.
\]
Therefore \eqref{eq:boundary} becomes
\[
 \partial_\eps e_i
 =\delta(g_i^+)[g_i^+]+\delta(g_i^-)[g_i^-].
\]
Let $\mathsf A$ be the unsigned incidence matrix, $\mathsf A e_i=[g_i^+]+[g_i^-]$, and let $D_\delta$ be the diagonal sign change on $V$ defined by $D_\delta[g]=\delta(g)[g]$. Since $D_\delta^2=\Id_V$, this map is invertible. Thus $\partial_\eps=D_\delta\mathsf A$ and $\dim\imageop\partial_\eps=\dim\imageop \mathsf A$.

A vector $y\in\ker \mathsf A^{\mathsf T}$ satisfies $y_u+y_v=0$ for each edge $uv$. On a connected bipartite component the values of $y$ alternate between the two vertex classes, giving one free parameter. An odd cycle forces that parameter to vanish. Thus $\dim\ker \mathsf A^{\mathsf T}=b(\Gamma_F)$ and $\dim\imageop \mathsf A=q-b(\Gamma_F)$, proving \eqref{eq:rank}. For constant orientations every edge joins germs of opposite sides, which supplies a bipartition of each component.
\end{proof}

For comparison, give each graph edge an arbitrary direction and assign to it its terminal vertex minus its initial vertex. The kernel of this ordinary incidence map is the \emph{graph cycle space}; its dimension is independent of the chosen directions. Write
\[
 \beta_1(\Gamma_F)=m-q+c(\Gamma_F),\qquad
 \nu(\Gamma_F)=c(\Gamma_F)-b(\Gamma_F)
\]
for this dimension and for the number of non-bipartite components, respectively.

\begin{corollary}[Parity defect of jump homology]\label{cor:parity}
For the jump graph,
\[
 r_\eps=q-c(\Gamma_F)+\nu(\Gamma_F),\qquad
 \dim H_1^\eps=\beta_1(\Gamma_F)-\nu(\Gamma_F).
\]
Each non-bipartite component contributes one to the difference between the ordinary and signed cycle dimensions.
\end{corollary}

\begin{corollary}[Complete three-branch jump-rank classification]\label{cor:three}
Suppose $m=2$ and both cuts are genuine jumps. For any nonconstant orientation pattern
\[
 --+ ,\quad -+- ,\quad -++ ,\quad +-- ,\quad +-+ ,\quad ++-,
\]
one has $r_\eps=2$ and $H_1^\eps=0$. For $---$ or $+++$, the rank is one exactly when the two cuts have the same left image germ and the same right image germ; otherwise the rank is two.
\end{corollary}
\begin{proof}
Each column of $\partial_\eps$ has exactly two nonzero entries, at distinct image germs. The possible pairs of germ sides for an edge at cut $i$ are $(-\eps_{i-1},\eps_i)$. For nonconstant three-branch patterns, the unordered side pairs at the two cuts differ. Thus their columns cannot be proportional. For constant orientations both side pairs consist of one positive and one negative germ; the two columns are proportional exactly when their corresponding left and right image germs coincide.
\end{proof}

For an alternating orientation pattern, let $\Gamma_F^\pm$ be the subgraphs on positive- and negative-side germs, with $q_\pm$ vertices and $b_\pm$ bipartite components.

\begin{corollary}[Alternating orientations]\label{cor:alternating}
If the orientations strictly alternate,
every edge joins germs of the same side, and
\[
 \Gamma_F=\Gamma_F^+\sqcup\Gamma_F^-,\qquad
 r_\eps=\bigl(q_+-b_+\bigr)+\bigl(q_--b_-\bigr).
\]
For at least three branches,
$r_\eps\geq2$.
\end{corollary}
\begin{proof}
At a cut where the orientation changes, both image germs have side
$\eps_i$. Under strict alternation no edge connects the positive-
and negative-side subgraphs, and with at least two cuts both subgraphs
have an edge. The rank assertions follow from Theorem~\ref{thm:rank}.
\end{proof}

\subsection{Integer coefficients and parity}

The boundary matrix has integer entries, so the chain complex of this section can also be formed over $\mathbb Z$. Its homology groups \cite[Section~2.1]{Hatcher2002} are
\[
 H_1^\varepsilon(\Gamma_F;\mathbb Z)
   =\ker(\partial_\varepsilon:\mathbb Z^m\to\mathbb Z^{\mathcal G}),
 \qquad
 H_0^\varepsilon(\Gamma_F;\mathbb Z)
   =\operatorname{coker}(\partial_\varepsilon:\mathbb Z^m\to\mathbb Z^{\mathcal G}).
\]
Here the cokernel is the quotient $\mathbb Z^{\calG}/\imageop\partial_\eps$. The same definitions apply over $\mathbb R$ or the two-element field $\mathbb F_2$; over $\mathbb R$ they recover $H_1^\eps$ and $H_0^\eps$ above. A summand $\mathbb Z/2\mathbb Z$ records an element of order two, invisible with real coefficients.

\begin{theorem}[Integral jump homology]\label{thm:integral-torsion}
For every piecewise strictly monotone interval map with finitely many
genuine jumps and arbitrary branch orientations,
\[
 H_0^\varepsilon(\Gamma_F;\mathbb Z)
 \cong \mathbb Z^{\,b(\Gamma_F)}
 \oplus (\mathbb Z/2\mathbb Z)^{\,\nu(\Gamma_F)},\qquad
 H_1^\varepsilon(\Gamma_F;\mathbb Z)
 \cong\mathbb Z^{\,m-|\mathcal G|+b(\Gamma_F)}.
\]
In particular,
\[
 \dim_{\mathbb F_2} H_1^\varepsilon(\Gamma_F;\mathbb F_2)
 -\dim_{\mathbb R}H_1^\varepsilon(\Gamma_F;\mathbb R)
 =\nu(\Gamma_F).
\]
\end{theorem}
\begin{proof}
Multiplying $\partial_\varepsilon$ by the diagonal matrix of
image-germ sides, which is invertible over $\mathbb Z$, on the left gives the
unsigned vertex--edge incidence matrix $\mathsf A$. The cokernel of $\mathsf A$ is
presented by generators $[v]$, one for each vertex, and relations
$[u]+[v]=0$ for each edge $uv$. On a connected component, choose a
spanning tree and a root $v_0$. The tree relations express every
vertex generator as $[v]=\pm[v_0]$. If the component is bipartite,
all remaining edge relations are consequences of the tree relations,
and its cokernel is $\mathbb Z$. Otherwise an odd cycle gives
$2[v_0]=0$, and every remaining edge relation is either trivial or
identical to this one. Its cokernel is thus $\mathbb Z/2\mathbb Z$.
Taking direct sums proves the assertion for $H_0^\varepsilon$.
The kernel of an integer matrix is free abelian, and its rank, computed after extension of scalars to $\mathbb Q$, is
$m-\dim_{\mathbb Q}\operatorname{im}\mathsf A=m-|\mathcal G|+b(\Gamma_F)$.
Finally, over $\mathbb F_2$, $\mathsf A$ is the ordinary oriented incidence
matrix, whose rank is $|\mathcal G|-c(\Gamma_F)$, whereas over
$\mathbb R$ its rank is $|\mathcal G|-b(\Gamma_F)$. Their kernel
dimensions differ by $c(\Gamma_F)-b(\Gamma_F)=\nu(\Gamma_F)$.
\end{proof}

\section{Cylinder growth and entropy}\label{sec:growth}

For $n\geq0$ and $1\leq i\leq m$, let $p_{n,i}$ count the points in $(a,b)$ whose first encounter with a cut occurs at time $n$, at $c_i$; thus $p_{0,i}=1$. Put
\[
 p_i(t)=\sum_{n\geq0}p_{n,i}t^n,\qquad
 p(t)=(p_1(t),\ldots,p_m(t))^{\mathsf T},\qquad \one=(1,\ldots,1)^{\mathsf T},
\]
and define the endpoint vector and cylinder series by
\[
 \xi(t)=\theta(b^-,t)-\theta(a^+,t),\qquad
 \mathscr L_F(t)=\sum_{n\geq0}\ell_n(F)t^n.
\]
The first identity below is the unweighted main kneading identity \cite[Theorem~2.1]{RughTan2015} in the present notation.

\begin{proposition}[Cylinder counting]\label{prop:counting}
For arbitrary branch orientations, the following formal identities hold:
\begin{align}
 N(t)p(t)&=\xi(t),\label{eq:mki}\\
 \mathscr L_F(t)&=\frac{1+t\one^{\mathsf T}N(t)^{-1}\xi(t)}{1-t},\label{eq:lap-resolvent}\\
 N(t)^{-1}&=\Id_E-tK_tR_\eps(t)^{-1}\partial_\eps.\label{eq:woodbury}
\end{align}
\end{proposition}
\begin{proof}
Fix a point $x$ whose first cut encounter is $F^n(x)=c_i$, and let $s\in\{-1,1\}$ be the orientation of $F^n$ on the surrounding $n$-cylinder. The first $n$ terms of the two one-sided series agree. Their tails have the common orientation factor $s$, while their image germs occur in the order $c_i^-,c_i^+$ if $s=1$ and in the reverse order if $s=-1$. Thus
\[
 \theta(x^+,t)-\theta(x^-,t)=t^n s^2N(t)e_i=t^nN(t)e_i.
\]
Each coefficient has only finitely many jumps. Summing them across $(a,b)$ proves \eqref{eq:mki} coefficientwise.

The points counted at time $k$ are the new boundaries of the $(k+1)$-cylinders, so
\[
 \ell_n(F)=1+\sum_{k=0}^{n-1}\sum_{i=1}^m p_{k,i}.
\]
Summing in $n$ gives \eqref{eq:lap-resolvent}. Multiplying the right-hand side of \eqref{eq:woodbury} by $\Id_E+tK_t\partial_\eps$, and using $R_\eps(t)=\Id_{\calB_\eps}+t\partial_\eps K_t|_{\calB_\eps}$, gives the identity operator.
\end{proof}

The radius of convergence of $\mathscr L_F$ is $\gamma_F^{-1}$. The next statement is the unweighted entropy theorem of Rugh and Tan, expressed through the reduced determinant.

\begin{theorem}[Entropy from the reduced determinant]\label{thm:entropy}\label{cor:entropy}
The function $D_F$ is analytic in $|t|<1$ and has no zero for $|t|<\gamma_F^{-1}$. If $\gamma_F>1$, then $\gamma_F^{-1}$ is its smallest positive zero. Consequently, $h_{\mathrm{top}}(F)=0$ if and only if $D_F$ has no positive zero in $(0,1)$.
\end{theorem}
\begin{proof}
The coefficients in \eqref{eq:theta} have absolute value at most $1/2$. Our matrix $N(t)$ is the transpose of the unweighted interior matrix $B(t)$ in \cite[Appendix~B]{RughTan2015}. Proposition~B.2 there gives
\[
 D_F(t)=(1-\alpha t)D_{\mathrm{RT}}(t),\qquad
 \alpha=\frac{\eps_0+\eps_m}{2},
\]
where $D_{\mathrm{RT}}$ is their full kneading determinant. Since $\alpha\in\{-1,0,1\}$, the factor $1-\alpha t$ has no zero in the open unit disc. Their Theorem~2.3 therefore gives the assertions.
\end{proof}

\section{Finite Markov systems}\label{sec:markov}

When the one-sided orbits of the cuts and endpoints lie in a finite set, interval transitions can be recorded by a finite matrix. The relation between its determinant and the kneading determinant includes factors from endpoint and germ dynamics. Related graph methods recover entropy for periodically driven Lorenz systems with finite critical orbits \cite{AlvesSilva2015}.

\subsection{Germ and interval chains}

Assume for the Markov identity that a finite ordered set
\[
 \calP=\{a=x_0<x_1<\cdots<x_N=b\}
\]
contains all cuts and is closed under all one-sided images. Such a set is a \emph{Markov partition}. Its partition intervals are $A_j=(x_j,x_{j+1})$ for $0\leq j<N$, and its \emph{transition matrix} is
\[
 M_{ij}=\begin{cases}
 1,&A_j\subset F(A_i),\\
 0,&\text{otherwise}.
 \end{cases}
\]
Each $A_i$ maps monotonically to an interval whose endpoints belong to $\calP$. Its contained partition intervals are recorded by row $i$ of $M$. Write $\rho(M)$ for the \emph{spectral radius}, the largest modulus of an eigenvalue of $M$.

Let $\widehat\calP$ consist of the two germs $x_j^\pm$ for every interior point and the existing endpoint germs $a^+$ and $b^-$. Set $W=\R^{\widehat\calP}$, with basis $[g]$ indexed by $g\in\widehat\calP$, and let $\partial_\eps:E\to W$ denote the boundary map \eqref{eq:boundary} followed by inclusion. Define
\[
 P_\eps[g]=\eps(g)[\ehat g],\qquad
 \Phi[g]=\varphi(g),\qquad
 T_\eps=P_\eps-\partial_\eps\Phi.
\]
Thus $P_\eps:W\to W$ follows signed germ images, $\Phi:W\to E$ records positions relative to the cuts, and $T_\eps:W\to W$ subtracts the jump contributions.

An element of $W$ is a \emph{germ chain}, a formal linear combination of the basis germs. Let $a_W:W\to\R$ be coefficient summation, $a_W([g])=1$. The subspace $W_0=\ker a_W$ consists of chains with total coefficient zero. For $1\leq j<N$, put
\[
 d_j=[x_j^+]-[x_j^-],\qquad Z=\spanop\{d_1,\ldots,d_{N-1}\}.
\]
Define $Q:Z\to Z$ by
\[
 Qd_j=\begin{cases}
 0,&x_j\text{ is a cut},\\
 d_k,&x_j\text{ is not a cut and }F(x_j)=x_k.
 \end{cases}
\]
The second case has $1\leq k<N$ by strict monotonicity. Let $u_0,\ldots,u_{N-1}$ be the standard basis of $\R^N$ and define $\eta:W\to\R^N$ by
\[
 \eta[x_j^\pm]=\sum_{k=0}^{j-1}u_k
\]
for existing germs at $x_j$, with the empty sum interpreted as zero. Thus $\eta$ replaces a germ by the interval chain from $a$ to its base point and forgets its side. In particular, $\eta d_j=0$. The chains
\[
 \zeta_j=[x_{j+1}^-]-[x_j^+]\quad(0\leq j<N)
\]
represent the partition intervals and satisfy $\eta\zeta_j=u_j$.

Put $\alpha=(\eps_0+\eps_m)/2$.
\begin{theorem}[Signed kneading--Markov identity]\label{thm:Markov} The subspaces $W_0$ and $Z$ are invariant under $T_\eps$, and
\[
 T_\eps|_Z=Q,\qquad
 \eta T_\eps|_{W_0}=M^{\mathsf T}\eta|_{W_0},\qquad
 W_0/Z\simeq\R^N.
\]
Moreover,
\begin{equation}\label{eq:Markov}
 D_F(t)=\frac{(1-\alpha t)\det(\Id_N-tM)
                    \det(\Id_Z-tQ)}
                   {\det(\Id_W-tP_\eps)}.
\end{equation}
Every zero of the two auxiliary determinants in \eqref{eq:Markov} lies on the unit circle. Consequently, for $|t|<1$, the zeros of $D_F(t)$ coincide, with multiplicity, with the zeros of $\det(\Id_N-tM)$, and
\[
 h_{\mathrm{top}}(F)=\log\rho(M).
\]
\end{theorem}
\begin{proof}
The signed germ orbit series gives $K_t=\Phi(\Id_W-tP_\eps)^{-1}|_V$. By Lemma~\ref{lem:recurrence} and Sylvester's identity $\det(\Id+AB)=\det(\Id+BA)$ for compatible rectangular matrices,
\begin{align*}
 D_F(t)&=\det\bigl(\Id_E+t\Phi(\Id_W-tP_\eps)^{-1}\partial_\eps\bigr)\\
 &=\frac{\det(\Id_W-tP_\eps+t\partial_\eps\Phi)}
         {\det(\Id_W-tP_\eps)}
 =\frac{\det(\Id_W-tT_\eps)}{\det(\Id_W-tP_\eps)}.
\end{align*}
Let $\chi_i(g)=\varphi_i(g)+1/2$ and let $\Phi_\chi[g]=(\chi_1(g),\ldots,\chi_m(g))^{\mathsf T}$. Since $\Phi=\Phi_\chi-\tfrac12\mathbf1 a_W$, with $\mathbf1=(1,\ldots,1)^{\mathsf T}\in E$, the operator $T_\eps$ satisfies
\begin{equation}\label{eq:uncentred}
 T_\eps=P_\eps-\partial_\eps\Phi_\chi
                +\tfrac12\partial_\eps\mathbf1\,a_W.
\end{equation}
For a germ in branch $j$, $a_WP_\eps[g]=\eps_j$. The telescoping identity
\[
 a_W\partial_\eps\Phi_\chi[g]
 =\sum_{i=1}^j(\eps_i-\eps_{i-1})=\eps_j-\eps_0
\]
and $a_W\partial_\eps\mathbf1=\eps_m-\eps_0$ give
\[
 a_W T_\eps=\alpha a_W.
\]
Hence $W_0$ is invariant and $T_\eps$ acts on $W/W_0$ by the scalar $\alpha$.

At a cut $x_j=c_i$, the difference of the two columns of $P_\eps$ is $\partial_\eps e_i$, while $\Phi d_j=e_i$; hence $T_\eps d_j=0$. At an interior non-cut point $\Phi d_j=0$. Germ transport gives $P_\eps d_j=d_k$ with $F(x_j)=x_k$: on a decreasing branch the orientation sign compensates for the interchange of the two sides. Therefore $Z$ is invariant and $T_\eps|_Z=Q$.

The sequence
\[
 0\longrightarrow Z\longrightarrow W_0
 \xrightarrow{\eta}\R^N\longrightarrow0
\]
is exact. Indeed, the $N-1$ chains $d_j$ and the $N$ chains $\zeta_j$ are the successive differences in the ordered list of $2N$ germs, so together they form a basis of $W_0$. The map $\eta$ kills the former and sends the latter to the basis $u_j$.

Both endpoints of $\zeta_j$ approach $A_j$ from within the same branch. Their centred indicators agree, so $\Phi\zeta_j=0$ and $T_\eps\zeta_j=P_\eps\zeta_j$. After applying $\eta$, the orientation sign makes the image chain run from the smaller image endpoint to the larger one. It is therefore the sum of the intervals contained in $F(A_j)$:
\[
 \eta T_\eps\zeta_j=\sum_{i=0}^{N-1}M_{ji}u_i
 =M^{\mathsf T}u_j.
\]
On each $d_j$, both $\eta T_\eps$ and $M^{\mathsf T}\eta$ vanish. The asserted intertwining identity follows on the whole basis of $W_0$.

The two invariant exact sequences give
\[
 \det(\Id_W-tT_\eps)=
 (1-\alpha t)\det(\Id_Z-tQ)\det(\Id_N-tM),
\]
proving \eqref{eq:Markov}. Each column of $P_\eps$ has a single entry equal to $+1$ or $-1$, and each column of $Q$ has a single entry $+1$ or is zero. Their nonzero eigenvalues therefore arise from finite cycles and are roots of unity. As $\alpha\in\{-1,0,1\}$, none of the auxiliary factors has a zero for $|t|<1$.

Finally, admissible words of length $n$ for the Markov partition are exactly the directed paths of length $n-1$ in $M$. Let $p_n$ count these paths and put $s=N-m-1$, the number of added interior partition points. Within an original $n$-cylinder, each $F^k$, $0\leq k<n$, is strictly monotone. Each added point therefore has at most one preimage under that iterate, so the refinement divides this cylinder into at most $1+ns$ intervals. Hence
\[
 \ell_n(F)\leq p_n\leq(1+ns)\ell_n(F).
\]
Path counts have exponential growth rate $\rho(M)$, and \eqref{eq:entropy-definition} gives the entropy assertion.
\end{proof}

\subsection{Itinerary-preserving parameter families}\label{sec:parameters}

The germ itineraries in Section~\ref{sec:setting} can be compared for
two maps with the same cuts and branch labels.
They record all successive one-sided branch choices, including those
after an orbit reaches a cut.

\begin{corollary}\label{cor:nonlinear}
Two maps satisfying the assumptions of Section~\ref{sec:setting},
with the same cuts, the same orientation word and the same itineraries of all one-sided cut
germs, have the same kneading determinant and the same topological
entropy.
Two maps satisfying the assumptions of Section~\ref{sec:setting}
also have the same entropy if they share a finite Markov partition
and the one-sided images of all its endpoints.
\end{corollary}
\begin{proof}
The first hypothesis makes each series $\theta(c_i^\pm,t)$ in
\eqref{eq:theta} identical for the two maps.
Theorem~\ref{thm:entropy} then gives equal entropy, including the
zero-entropy case.
In the second assertion the transition matrices agree, so
Theorem~\ref{thm:Markov} applies.
\end{proof}

Let $\Lambda\subset\R^d$, with $d\ge1$, and let
$(F_\lambda)_{\lambda\in\Lambda}$ be a family satisfying the standing
assumptions. Assume that its orientation word is fixed. Its \emph{incidence pattern} records which of the labelled
germs $g_i^\pm$ coincide. Where this pattern is fixed, the graph
spaces $E,V,\calB_\eps$ are naturally identified and the reduced determinant
has fixed order $r_\eps$.

\section{Orbital symmetry}\label{sec:symmetry}

\subsection{Compatible involutions}

An involution is a linear operator whose square is the identity. Let $J:E\to E$ and $S_\calB:\calB_\eps\to\calB_\eps$ be involutions satisfying
\begin{equation}\label{eq:symmetry-compatibility}
 \partial_\eps J=S_\calB\partial_\eps,\qquad
 K_tS_\calB=JK_t\quad\hbox{on }\calB_\eps.
\end{equation}
We call them compatible with the orbital data, and write
\[
 E^\pm=\ker(J\mp\Id_E),\qquad
 \calB_\eps^\pm=\ker(S_\calB\mp\Id_{\calB_\eps}).
\]

Compatibility gives $N(t)J=JN(t)$ and $R_\eps(t)S_\calB=S_\calB R_\eps(t)$, so the eigenspaces are invariant. Define
\[
 D_F^\pm(t)=\det\bigl(R_\eps(t)|_{\calB_\eps^\pm}\bigr).
\]
\begin{proposition}[Orbital factorisation]\label{prop:symmetry}
The sector determinants satisfy
\begin{equation}\label{eq:sym}
 D_F^\pm(t)=\det\bigl(N(t)|_{E^\pm}\bigr),\qquad
 D_F(t)=D_F^+(t)D_F^-(t).
\end{equation}
\end{proposition}
\begin{proof}
The boundary sequence \eqref{eq:boundary-exact} restricts to each eigenspace: projecting a lift by $(\Id_E\pm J)/2$ proves surjectivity onto $\calB_\eps^\pm$. On each kernel the kneading operator is the identity. Taking determinants proves \eqref{eq:sym}.
\end{proof}

\begin{corollary}[Entropy in an invariant sector]\label{cor:symmetric-entropy}
Assume \eqref{eq:symmetry-compatibility}, $J\xi(t)=\xi(t)$ and $\one^{\mathsf T}J=\one^{\mathsf T}$. If $h_{\mathrm{top}}(F)>0$, then $\gamma_F^{-1}$ is the smallest positive zero of $D_F^+(t)$.
\end{corollary}
\begin{proof}
Since $N(t)$ commutes with $J$, the series $N(t)^{-1}\xi(t)$ belongs to $E^+[[t]]$. Formula~\eqref{eq:lap-resolvent} therefore depends only on the inverse of $N(t)|_{E^+}$. If $D_F^+(\gamma_F^{-1})\ne0$, the restricted inverse would be analytic there, and that formula would analytically continue $\mathscr L_F$ through its positive radius of convergence. This contradicts Pringsheim's theorem \cite[Theorem~IV.6]{FlajoletSedgewick2009}. Theorem~\ref{thm:entropy} excludes an earlier zero, and Proposition~\ref{prop:symmetry} identifies the boundary factor.
\end{proof}

\subsection{Reflection}

Assume $I=[-b,b]$, $c_i=-c_{m+1-i}$, and $F(-x)=-F(x)$ on its branch domain. Then $\eps_j=\eps_{m-j}$. Reflection acts on germs by $\mathfrak r(x^\sigma)=(-x)^{-\sigma}$; it preserves $\calG$. Put $R_V[g]=[\mathfrak r g]$ and $Je_i=e_{m+1-i}$.

\begin{lemma}[Reflection compatibility]\label{lem:reflection}
The boundary and itinerary operators satisfy
\[
 R_V\partial_\eps=-\partial_\eps J,\qquad
 \theta(\mathfrak r g,t)=-J\theta(g,t),\qquad K_tR_V=-JK_t.
\]
\end{lemma}
\begin{proof}
Oddness gives $\ehat\mathfrak r=\mathfrak r\ehat$ and $\eps(\mathfrak r g)=\eps(g)$, while $\varphi(\mathfrak r g)=-J\varphi(g)$. Summing \eqref{eq:theta} gives the itinerary identity. The boundary identity follows from $\mathfrak r g_i^+=g_{m+1-i}^-$, $\mathfrak r g_i^-=g_{m+1-i}^+$ and \eqref{eq:boundary}.
\end{proof}

Thus $S_\calB=-R_V|_{\calB_\eps}$ satisfies \eqref{eq:symmetry-compatibility}. The terms symmetric and antisymmetric refer to the $+$ and $-$ eigenspaces of the cut reversal $J$ and their boundary images. The determinant on a zero-dimensional space is one.

\begin{corollary}[Reflection and entropy]\label{thm:sym}\label{cor:leading}
For an odd map with symmetric cuts, $D_F=D_F^+D_F^-$. If $h_{\mathrm{top}}(F)>0$, its smallest positive zero belongs to $D_F^+$. For two symmetric cuts, both factors have order at most one; if $r_\eps=1$, then $D_F^-=1$.
\end{corollary}
\begin{proof}
Reflection interchanges $a^+$ and $b^-$, so Lemma~\ref{lem:reflection} gives $J\xi=\xi$. Also $\one^{\mathsf T}J=\one^{\mathsf T}$. Apply Proposition~\ref{prop:symmetry} and Corollary~\ref{cor:symmetric-entropy}. When $m=2$, the two cut eigenspaces are spanned by $e_1+e_2$ and $e_1-e_2$. If $r_\eps=1$, Corollary~\ref{cor:three} gives $\partial_\eps(e_1-e_2)=0$, so the antisymmetric factor is one.
\end{proof}

Global oddness can be replaced by conditions on the relevant orbit data: it suffices that reflection preserve $\calG$, that the boundary identity hold, and that $\theta(\mathfrak r g,t)=-J\theta(g,t)$ hold for the jump-image and endpoint germs. These conditions imply the same compatibility and endpoint invariance even if $F$ is not odd elsewhere.

For a finite symmetric Markov partition, define $J_Mu_j=u_{N-1-j}$, $R_W[g]=[\mathfrak r g]$ and $W^-=\ker(R_W+\Id_W)$. Reflection commutes with $P_\eps$ and $T_\eps$, while $\Phi R_W=-J\Phi$. Hence $\partial_\eps(E^+)\subset W^-$ and $\Phi(W^-)\subset E^+$. The determinant calculation in Theorem~\ref{thm:Markov}, restricted to these spaces, gives
\begin{equation}\label{eq:symmetric-Markov}
 D_F^+(t)=\frac{\det(\Id_{W^-}-tT_\eps|_{W^-})}
                    {\det(\Id_{W^-}-tP_\eps|_{W^-})}.
\end{equation}
The denominator has no zero in $|t|<1$. In particular, when $\rho(M)>1$, the positive zero in this factor is $\rho(M)^{-1}$.

\subsection{Cumulative orientation}\label{subsec:cover}

For a Markov partition as in Section~\ref{sec:markov}, let $M_+$ be the $N\times N$ matrix obtained from $M$ by retaining the rows arising from increasing branches and setting the other rows to zero; let $M_-=M-M_+$. To record cumulative orientation, replace each state $A_j$ by $(A_j,+)$ and $(A_j,-)$. A transition preserves this extra sign on an increasing branch and reverses it on a decreasing branch. The resulting \emph{orientation double cover} has transition matrix $\mathcal M$, where
\[
 S=M_+-M_-,\qquad
 \mathcal M=\begin{pmatrix}M_+&M_-\\M_-&M_+\end{pmatrix}.
\]
Here $S$ is the transition matrix weighted by branch orientation. Let $U(u,v)=(u+v,u-v)$ on $\R^N\oplus\R^N$; its inverse sends $(p,q)$ to $((p+q)/2,(p-q)/2)$.

\begin{proposition}[Orientation decomposition]\label{prop:cover}
The map $U$ conjugates $\mathcal M$ to $\operatorname{diag}(M,S)$. Thus
\begin{equation}\label{eq:double}
 \det(\Id_{2N}-t\mathcal M)
 =\det(\Id_N-tM)\det(\Id_N-tS).
\end{equation}
If every branch is decreasing, $\mathcal M^2=\operatorname{diag}(M^2,M^2)$; if every branch is increasing, $\mathcal M=\operatorname{diag}(M,M)$.
\end{proposition}
\begin{proof}
On a vector $(u,u)$, $\mathcal M$ acts by $M u$ in both components. On $(u,-u)$, it acts by $(Su,-Su)$. These two invariant subspaces are complementary and give \eqref{eq:double}. The special cases follow by setting $(M_+,M_-)=(0,M)$ or $(M,0)$.
\end{proof}

For an odd Markov system with a symmetric partition, put $E_M^\pm=\ker(J_M\mp\Id_N)$.

\begin{corollary}[Four symmetry sectors]\label{cor:four}
Both $M$ and $S$ preserve $E_M^+$ and $E_M^-$. The orientation cover therefore has the determinant factorisation
\[
 \det(\Id_{2N}-t\mathcal M)=
 \prod_{\sigma\in\{+,-\}}
 \det(\Id_{E_M^\sigma}-tM|_{E_M^\sigma})
 \det(\Id_{E_M^\sigma}-tS|_{E_M^\sigma}).
\]
\end{corollary}
\begin{proof}
Reflection interchanges symmetric partition intervals, and oddness
preserves their orientation signs. Hence $J_MM_\pm=M_\pm J_M$.
Both $M=M_++M_-$ and $S=M_+-M_-$ commute with $J_M$.
Apply Proposition~\ref{prop:cover} on the two eigenspaces of $J_M$.
\end{proof}

\section{Sharp rank and entropy realisations}\label{sec:sharp}

\subsection{The smallest rank for an orientation word}

Fix an orientation word $\varepsilon=(\varepsilon_0,\ldots,\varepsilon_m)$. Let $\mathcal F_\varepsilon$ be the class of maps satisfying the assumptions of Section~\ref{sec:setting}, with this orientation word and all cuts genuine jumps; both the interval $I$ and the positions of the cuts may vary. Throughout this subsection and the next, $n=m+1$ denotes the number of branches.
At $c_i$, the adjacent branch-orientation pair $(\varepsilon_{i-1},\varepsilon_i)$ gives the image-germ side pair $(\delta(g_i^-),\delta(g_i^+))=(-\varepsilon_{i-1},\varepsilon_i)$. We classify the branch-orientation pairs into three \emph{side types} according to these image-germ sides:
\[
\begin{aligned}
 \mathsf P&:\ (-,+)\quad\text{(two positive germs)},\\
 \mathsf N&:\ (+,-)\quad\text{(two negative germs)},\\
 \mathsf X&:\ (--),(++)\quad\text{(one germ of each side)}.
\end{aligned}
\]
Let $d(\varepsilon)$ be the number of these types represented among
the $m$ cuts, so $1\le d(\varepsilon)\le3$.

\begin{theorem}[Sharp minimum of the homological rank]
\label{thm:sharp-pattern-rank}
For every orientation word $\varepsilon\in\{-1,+1\}^{m+1}$,
\[
 \min_{F\in\mathcal F_\varepsilon} r_\varepsilon(F)=d(\varepsilon).
\]
The minimum is attained by a piecewise affine Markov map on $[0,1]$.
\end{theorem}
\begin{proof}
Multiplication by $D_\delta$ changes the signed boundary into unsigned incidence columns. A column of type $\mathsf P$ has two distinct positive-side vertices, a column of type $\mathsf N$ has two distinct negative-side vertices, and a column of type $\mathsf X$ has one vertex of each side. Distinctness follows from the genuine-jump assumption.

Choose one column of each represented type. Any two of different types are linearly independent: their supports either lie on opposite sides or have different numbers of vertices on each side. If all three types occur, consider a relation among their chosen columns. Its positive-side projection is a relation between the $\mathsf P$ column, with two nonzero coordinates, and the positive part of the $\mathsf X$ column, with only one. These two vectors are independent, so both coefficients vanish. The remaining coefficient of the $\mathsf N$ column also vanishes. This proves the lower bound in all cases.

For sharpness put $c_i=i/n$, $I=[0,1]$ and prescribe the
one-sided images of the cuts as follows:
\[
\begin{array}{c|cc}
(\varepsilon_{i-1},\varepsilon_i)&F(c_i^-)&F(c_i^+)\\ \hline
(-,-)&0&1\\
(+,+)&1&0\\
(-,+)&0&c_1\\
(+,-)&1&c_m.
\end{array}
\]
Set $F(0^+)=1$ if $\varepsilon_0=-1$ and $F(0^+)=0$ otherwise;
set $F(1^-)=0$ if $\varepsilon_m=-1$ and $F(1^-)=1$ otherwise.
Interpolate affinely within each open branch. Every increasing branch
has left image at most $c_1$ and right image $1$, whereas every
decreasing branch has left image at least $c_m$ and right image $0$.
Thus the prescribed branches are strictly monotone in the required
directions. All jumps are genuine.

Every edge of type $\mathsf P$ joins the germ vertices $0^+$ and $c_1^+$,
every edge of type $\mathsf N$ joins $1^-$ and $c_m^-$, and every edge of type
$\mathsf X$ joins $0^+$ and $1^-$. After parallel edges are identified, the
jump graph is a subgraph of the path
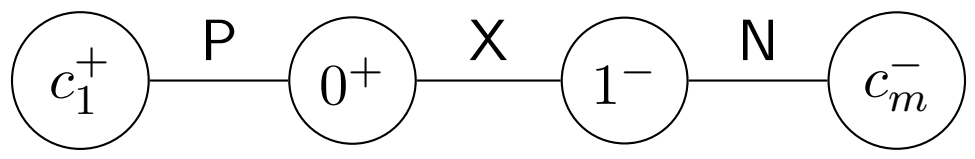
\[
\begin{tikzpicture}[baseline=-2pt]
 \foreach \name/\x/\lab in {A/0/{c_1^+},B/1.8/{0^+},C/3.6/{1^-},D/5.4/{c_m^-}}
   \node[circle,draw,inner sep=3pt] (\name) at (\x,0) {$\lab$};
 \draw (A)--node[above] {$\mathsf P$}(B)--node[above] {$\mathsf X$}(C)--node[above] {$\mathsf N$}(D);
\end{tikzpicture}
\]
It is an acyclic graph, or forest, with exactly $d(\varepsilon)$ distinct edges, hence
$\dim\operatorname{im}\partial_\varepsilon=d(\varepsilon)$.
Finally, all one-sided images of $\{0,c_1,\ldots,c_m,1\}$ belong to
this same set, so it is a Markov partition.
\end{proof}

\subsection{Entropy of the realisations}

The affine realisation also permits an explicit entropy calculation.
Write $A$ for the number of transitions $-\to+$, $B$ for the number
of transitions $+\to-$, and $\kappa=A+B$. The final adjacent pair determines
\[
 a_*=\begin{cases}1,&(\eps_{m-1},\eps_m)=(-,+),\\0,&\text{otherwise},\end{cases}
 \qquad
 b_*=\begin{cases}1,&(\eps_{m-1},\eps_m)=(+,-),\\0,&\text{otherwise}.\end{cases}
\]

\begin{theorem}[Uniform cubic Markov reduction and high entropy]
\label{thm:sharp-high-entropy}
For the affine Markov map constructed in
Theorem~\ref{thm:sharp-pattern-rank}, the Markov matrix $M$ has
size $n$ and rank at most three. Its Markov determinant
satisfies
\[
 \det(\Id_n-tM)=
 1-(n-b_*)t+(\kappa-nb_*)t^2+(Ab_*-Ba_*)t^3.
\]
For every $n\ge2$,
\[
 0\leq\log n-h_{\mathrm{top}}(F)
 \leq\log\frac{n}{n-1}\longrightarrow0\qquad(n\to\infty),
\]
uniformly over all orientation words.
Consequently, any prescribed family of orientation words of increasing
length admits maps with minimal homological rank at most three and
unbounded topological entropy.
\end{theorem}
\begin{proof}
Every branch image contains the interior partition intervals indexed by $1,\ldots,n-2$, so the corresponding entries in every row of $M$ are one. The only possible missing
entry in column $0$ occurs in a row $j$ where the orientation changes
from $-$ to $+$ at $c_j$; the only possible missing entry in column
$n-1$ occurs where the orientation changes from $+$ to $-$.
Let $u,v\in\{0,1\}^n$ be their indicator columns, let $e_0,e_{n-1}$
be standard coordinate vectors in $\R^n$, and write
$\mathbf1=(1,\ldots,1)^{\mathsf T}\in\R^n$ in this proof. Then
\[
 M=\mathbf1\mathbf1^{\mathsf T}-u e_0^{\mathsf T}
                     -v e_{n-1}^{\mathsf T}.
\]
Here $\mathbf1^{\mathsf T}u=A$, $\mathbf1^{\mathsf T}v=B$,
$u_0=v_0=0$, $u_{n-1}=a_*$ and $v_{n-1}=b_*$. Thus $M=X_0Y_0$ with
\[
 X_0=(\mathbf1,u,v),\qquad
 Y_0=\begin{pmatrix}
 \mathbf1^{\mathsf T}\\-e_0^{\mathsf T}\\-e_{n-1}^{\mathsf T}
 \end{pmatrix},\qquad
 Y_0X_0=\begin{pmatrix}
 n&A&B\\-1&0&0\\-1&-a_*&-b_*
 \end{pmatrix}.
\]
Sylvester's determinant identity gives
$\det(\Id_n-tM)=\det(\Id_3-tY_0X_0)$; expansion proves the displayed cubic.
The supports of $u$ and $v$ are disjoint, since a cut has only one orientation-transition type. Each row sum of $M$ is therefore $n$ or $n-1$. The spectral-radius bounds by minimum and maximum row sums give
\[
 n-1\leq\rho(M)\leq n.
\]
The Markov entropy formula completes the proof.
\end{proof}

\begin{corollary}[Entropy under rearrangement]
\label{cor:entropic-universality}
If the final two orientations are equal, then $a_*=b_*=0$ and
\[
 \det(\Id_n-tM)=1-nt+\kappa t^2,\qquad
 h_{\mathrm{top}}(F)
 =\log\left(\frac{n+\sqrt{n^2-4\kappa}}{2}\right).
\]
The entropy depends only on $n$ and $\kappa$ within this subfamily.
\end{corollary}
\begin{proof}
Setting $a_*=b_*=0$ in Theorem~\ref{thm:sharp-high-entropy} gives the determinant. Here $0\leq\kappa\leq n-2$, so $n^2-4\kappa>0$. The larger root of $\lambda^2-n\lambda+\kappa=0$ is the Perron eigenvalue and gives the entropy formula.
\end{proof}

For example, let $\varepsilon^{(1)}=(-,+,+,-,-,-)$ and $\varepsilon^{(2)}=(+,+,-,-,+,+)$. Write $M_1,M_2$ for their Markov matrices and
\[
 S_j=\operatorname{diag}(\varepsilon^{(j)})M_j\qquad(j=1,2)
\]
for their orientation-weighted matrices, as in Section~\ref{subsec:cover}. Both patterns have $\kappa=2$, so
\[
 \det(\Id_6-tM_1)=\det(\Id_6-tM_2)=1-6t+2t^2,\qquad
 h_{\mathrm{top}}(F)=\log(3+\sqrt7)
\]
for either map. Applying the three-column factorisation to $S_1$ and $S_2$ gives the reduced matrices
\[
 \begin{pmatrix}-2&1&1\\1&0&0\\-1&0&0\end{pmatrix}
 \quad\text{and}\quad
 \begin{pmatrix}2&1&1\\-1&0&0\\1&0&0\end{pmatrix}.
\]
Consequently,
\[
 \det(\Id_6-tS_1)=1+2t,\qquad \det(\Id_6-tS_2)=1-2t,
\]
and, for every integer $p\geq1$,
\[
 \operatorname{tr}(S_1^p)=(-2)^p,\qquad
 \operatorname{tr}(S_2^p)=2^p.
\]
The traces count closed Markov paths with their cumulative orientation signs. Thus these maps have identical nonzero Markov spectra and different signed counts for every odd period.

\begin{proposition}[Maximal entropy at minimal rank]
\label{prop:full-entropy-minimal-rank}
For every orientation word $\varepsilon\in\{-1,+1\}^{n}$, there exists
$F\in\mathcal F_\varepsilon$ such that
\[
 r_\varepsilon(F)=d(\varepsilon),\qquad
 \ell_k(F)=n^k\quad(k\ge0).
\]
Consequently, $h_{\mathrm{top}}(F)=\log n$ and
$\mathscr L_F(t)=(1-nt)^{-1}$.
\end{proposition}
\begin{proof}
Choose $0<\delta<\eta<1/n$. Keep the cuts and endpoint values from the
construction in Theorem~\ref{thm:sharp-pattern-rank}, replacing the
cut-image values $c_1,c_m$ in its table by $\delta,1-\delta$, respectively.
On the outer branches also prescribe
\[
 F(\eta)=\frac{1-\varepsilon_0}{2}+\varepsilon_0\delta,
 \qquad
 F(1-\eta)=\frac{1+\varepsilon_m}{2}-\varepsilon_m\delta.
\]
Interpolate affinely between successive prescribed points on each branch.
The branches remain strictly monotone with orientations $\varepsilon$,
and all jumps are genuine. The jump graph retains the same incidence
pattern, so $r_\varepsilon(F)=d(\varepsilon)$.

Put $K=[\eta,1-\eta]$. The image of each $K\cap J_j$ has one endpoint
at most $\delta$ and the other at least $1-\delta$. Hence
$K\subset F(K\cap J_j)$ for every $0\le j<n$.
Each inverse branch therefore maps $K$ into $K\cap J_j$.
Composing these inverse branches realizes every finite branch word,
so $\ell_k(F)=n^k$. The entropy and cylinder series follow.
\end{proof}

\section{Applications with mixed orientations}\label{sec:mixed}

\subsection{A symmetric three-branch map with orientation \texorpdfstring{$(-,+,-)$}{(-,+,-)}}

Let $I=[0,1]$, with cuts $1/3$ and $2/3$, and define
\begin{equation}\label{eq:example}
 F(x)=\begin{cases}
 1-2x,&0<x<1/3,\\
 3x-1,&1/3<x<2/3,\\
 2-2x,&2/3<x<1.
 \end{cases}
\end{equation}
This system satisfies $F(1-x)=1-F(x)$, has two genuine jumps, and has orientation $(-,+,-)$. In centred coordinates, $\widetilde F(z)=F(z+1/2)-1/2$ is odd on $[-1/2,1/2]$, so Section~\ref{sec:symmetry} applies. Its cut-image germs are
\[
 g_1^-=(1/3)^+ ,\quad g_1^+=0^+ ,\qquad
 g_2^-=1^- ,\quad g_2^+=(2/3)^- .
\]
The signed jump graph has two disjoint edges, so $r_\eps=2$. The original branch partition is Markov, with
\[
 M=\begin{pmatrix}0&1&1\\1&1&1\\1&1&0\end{pmatrix},\qquad
 \det(\Id_3-tM)=(1+t)(1-2t-t^2),\qquad
 \rho(M)=1+\sqrt2.
\]
The endpoint germs form the cycle $0^+\mapsto1^-\mapsto0^+$, with orientation $-1$ at both steps. The other two image germs enter this cycle after one step. Applying \eqref{eq:recurrence} first on the cycle and then at the cuts gives
\[
 N(t)=\frac{1}{1-t}
 \begin{pmatrix}1-t-t^2&-t\\-t&1-t-t^2\end{pmatrix}.
\]
Therefore
\begin{equation}\label{eq:exampleD}
 D_F(t)=\frac{(1+t)(1-2t-t^2)}{1-t},\qquad
 D_F^+(t)=\frac{1-2t-t^2}{1-t},\qquad
 D_F^-(t)=1+t,
\end{equation}
so
\[
 h_{\mathrm{top}}(F)=\log(1+\sqrt2).
\]
The same result follows from Theorem~\ref{thm:Markov}, because $\alpha=-1$, $Q=0$, and $\det(\Id_W-tP_\eps)=1-t^2$. Reflection splits the two-dimensional determinant into two scalar factors.

\subsection{A four-branch alternating reduction}

For the orientation pattern $(-,+,-,+)$ on $[0,1]$, take cuts
$c_1=1/4$, $c_2=1/2$, $c_3=3/4$ and the unique affine branches
with endpoint values
\[
 \begin{array}{c|cc}
 \text{branch}&\text{left endpoint image}&\text{right endpoint image}\\\hline
 0&9/10&1/5\\
 1&2/5&4/5\\
 2&7/10&1/5\\
 3&2/5&9/10
 \end{array}
\]
The first and third jumps have identical positive-side image germs,
while the second jump has two distinct negative-side image germs.
Consequently $\partial_\eps e_1=\partial_\eps e_3$, whereas
$\partial_\eps e_2$ has disjoint support. Thus
\[
 \ker\partial_\eps=\spanop\{e_1-e_3\},\qquad r_\eps=2<3.
\]
The alternating pattern can therefore exhibit genuine homological
reduction once there are at least four branches. Its jump graph consists of two parallel positive-side edges and one negative-side edge.

\subsection{An odd cycle in the jump graph}

Consider four affine branches on $[0,1]$, with cuts $1/4$, $1/2$, and $3/4$, having orientations $(-,+,+,+)$ and endpoint values
\[
 \begin{array}{c|cc}
 \text{branch}&\text{left endpoint image}&\text{right endpoint image}\\\hline
 0&9/10&1/5\\
 1&3/10&4/5\\
 2&3/10&4/5\\
 3&1/5&9/10
 \end{array}
\]
(the indicated values are one-sided at cuts). Each branch has a unique affine extension. The three jump edges have image-germ endpoints
\[
 (1/5)^+\longleftrightarrow(3/10)^+,
 \qquad (4/5)^-\longleftrightarrow(3/10)^+,
 \qquad (4/5)^-\longleftrightarrow(1/5)^+.
\]
The jump graph is a triangle. Thus $q=3$, $b(\Gamma_F)=0$, and $r_\eps=3$, despite having only three distinct jump-image germs. This illustrates why the bipartite-component count, rather than the ordinary connected-component count, governs the general homological rank.

\section{Applications to decreasing branches}\label{sec:decreasing}\label{sec:common}

Assume throughout this section that every branch is decreasing. Put $\partial=-\partial_\eps$, so $\partial e_i=[g_i^+]-[g_i^-]$ and $R_\eps(t)=\Id_{\calB_\eps}-t\partial K_t|_{\calB_\eps}$. The \emph{uncentred indicator series} is
\[
 \theta^{\mathrm u}(g,t)
 =\sum_{n\ge0}(-t)^n\bigl(\varphi(\ehat^n g)+\tfrac12\one\bigr).
\]
Thus $\theta^{\mathrm u}(g,t)=\theta(g,t)+\one/[2(1+t)]$ for every germ. The added vector cancels in all germ differences and on $\calB_\eps$.

\subsection{Common endpoint limits}

Assume that every jump has the same endpoint limits:
\begin{equation}\label{eq:full-resets}
 F(c_i^-)=a,\qquad F(c_i^+)=b,\qquad 1\le i\le m.
\end{equation}
The jump graph has two vertices and $m$ parallel edges.
Thus $r_\eps=1$ and $\dim H_1^\eps=m-1$.
Put $\Delta=[b^-]-[a^+]$. Then $\partial e_i=\Delta$ and
$K_t\Delta=\xi(t)$.
The branch address satisfies $\iota(g)=\sum_{i=1}^m(\varphi_i(g)+1/2)$.
Define the scalar endpoint series by
\begin{equation}\label{eq:scalar-series}
 K_F(t)=\one^{\mathsf T}\xi(t)
 =\sum_{n\ge0}(-t)^n
   \bigl(\iota(\widehat F^nb^-)-\iota(\widehat F^na^+)\bigr).
\end{equation}

\begin{corollary}\label{cor:rank-one}
Under \eqref{eq:full-resets},
\begin{equation}\label{eq:scalar-formula}
 D_F(t)=1-tK_F(t),\qquad
 \mathscr L_F(t)=\frac{1}{(1-t)D_F(t)}.
\end{equation}
If $m\ge3$, then
$\log(m-1)\le h_{\mathrm{top}}(F)\le\log(m+1)$ and
$\gamma_F^{-1}\in[1/(m+1),1/(m-1)]$.
\end{corollary}
\begin{proof}
Here $N(t)=\Id_E-t\,\xi(t)\one^{\mathsf T}$, so
$\det N(t)=1-t\one^{\mathsf T}\xi(t)$ and
$N(t)^{-1}\xi=\xi/D_F$.
Substitution in \eqref{eq:lap-resolvent} proves
\eqref{eq:scalar-formula}.
Every word on the $m-1$ interior branches occurs, since each such
branch maps onto $(a,b)$. There are at most $(m+1)^n$ words of
length $n$. The entropy bounds follow from
\eqref{eq:entropy-definition}, and the zero bounds follow from
Theorem~\ref{thm:entropy}.
\end{proof}

\subsection{Odd bimodal maps}

Under \eqref{eq:full-resets}, suppose $m=2$, $a=-b$, the two cuts are symmetric, and $F$ is odd,
meaning $F(-x)=-F(x)$ away from the cuts.
This is the odd discontinuous bimodal setting of \cite{Oliveira2015}.
Reflection gives
$\iota(\widehat F^na^+)=2-\iota(\widehat F^nb^-)$, and hence
\begin{equation}\label{eq:2015-scalar}
 D_F(t)=1-2t\sum_{n\ge0}(-t)^n
             \bigl(\iota(\widehat F^nb^-)-1\bigr).
\end{equation}
The normalisation used in \cite{Oliveira2015} is $D_F(t)/(1+t)$.
The one-sequence formula thus follows from reflection of the endpoint
orbits, without a periodicity assumption.
When these orbits are finite, Theorem~\ref{thm:Markov}
also gives the characteristic-polynomial relation.

For the motivating tangent family $f_\beta(x)=-\beta\tanh(\beta\tan x)$, $\pi/2<\beta<3\pi/2$, both cuts $\pm\pi/2$ have left limit $-\beta$ and right limit $\beta$. The derivative is negative on every open branch and the map is odd. Formula~\eqref{eq:2015-scalar} therefore applies directly.

\subsection{Four equal jumps and two independent parameters}\label{sec:decreasing-examples}
\begin{figure}[htbp]
\centering
\begin{minipage}[c]{.49\textwidth}
\centering
\begin{tikzpicture}[x=2.35cm,y=2.35cm,font=\small]
 \draw[map axis] (-1,0)--(1,0) (0,-1)--(0,1);
 \draw[map frame] (-1,-1) rectangle (1,1);
 \foreach \x/\lab in {-1/{-1},-.75/{-\frac34},-.25/{-\frac14},
                       .25/{\frac14},.75/{\frac34},1/{1}}{
  \draw[map frame] (\x,-1)--(\x,-1.035);
  \node[below,inner sep=2pt] at (\x,-1.035) {$\lab$};
 }
 \foreach \y in {-1,0,1}{
  \draw[map frame] (-1,\y)--(-1.035,\y);
  \node[left,inner sep=2pt] at (-1.035,\y) {$\y$};
 }
 \node[anchor=south west] at (-1,1.07) {$F(x)$};
 \node[anchor=west] at (1.055,-1) {$x$};
 \node[anchor=south east] at (1,1.07) {(a)};
 \draw[map branch] (-1,.25)--(-.75,-1);
 \draw[map branch] (-.75,1)--(-.25,-1);
 \draw[map branch] (-.25,1)--(.25,-1);
 \draw[map branch] (.25,1)--(.75,-1);
 \draw[map branch] (.75,1)--(1,-.75);
 \foreach \x in {-.75,-.25,.25,.75}{
  \node[open end] at (\x,-1) {};
  \node[open end] at (\x,1) {};
 }
 \node[closed end] at (-1,.25) {};
 \node[closed end] at (1,-.75) {};
\end{tikzpicture}
\end{minipage}\hfill
\begin{minipage}[c]{.49\textwidth}
\centering
\begin{tikzpicture}[font=\small]
 \path[use as bounding box] (-.6,-2.8) rectangle (4.4,2.8);
 \node[anchor=north east] at (4.3,2.7) {(b)};
 \coordinate (a) at (0,0);
 \coordinate (b) at (3.8,0);
 \draw[jump edge] (a) .. controls (.8,2.1) and (3,2.1) ..
  node[edge label,pos=.5] {$e_1$} (b);
 \draw[jump edge] (a) .. controls (1,.7) and (2.8,.7) ..
  node[edge label,pos=.5] {$e_2$} (b);
 \draw[jump edge] (a) .. controls (1,-.7) and (2.8,-.7) ..
  node[edge label,pos=.5] {$e_3$} (b);
 \draw[jump edge] (a) .. controls (.8,-2.1) and (3,-2.1) ..
  node[edge label,pos=.5] {$e_4$} (b);
 \node[closed end] at (a) {};
 \node[closed end] at (b) {};
 \node[left=4pt] at (a) {$(-1)^+$};
 \node[right=4pt] at (b) {$1^-$};
 \node at (1.9,-2.25) {$\dim\calB_\eps=1,\qquad \dim H_1^\eps=3$};
\end{tikzpicture}
\end{minipage}
\caption{(a) The non-odd map $F_{1/4,-3/4}$.
(b) Its four jump edges have the same initial and terminal germs,
so the reduced determinant is scalar.
Open circles denote excluded one-sided branch endpoints.}
\label{fig:rank-one}
\end{figure}
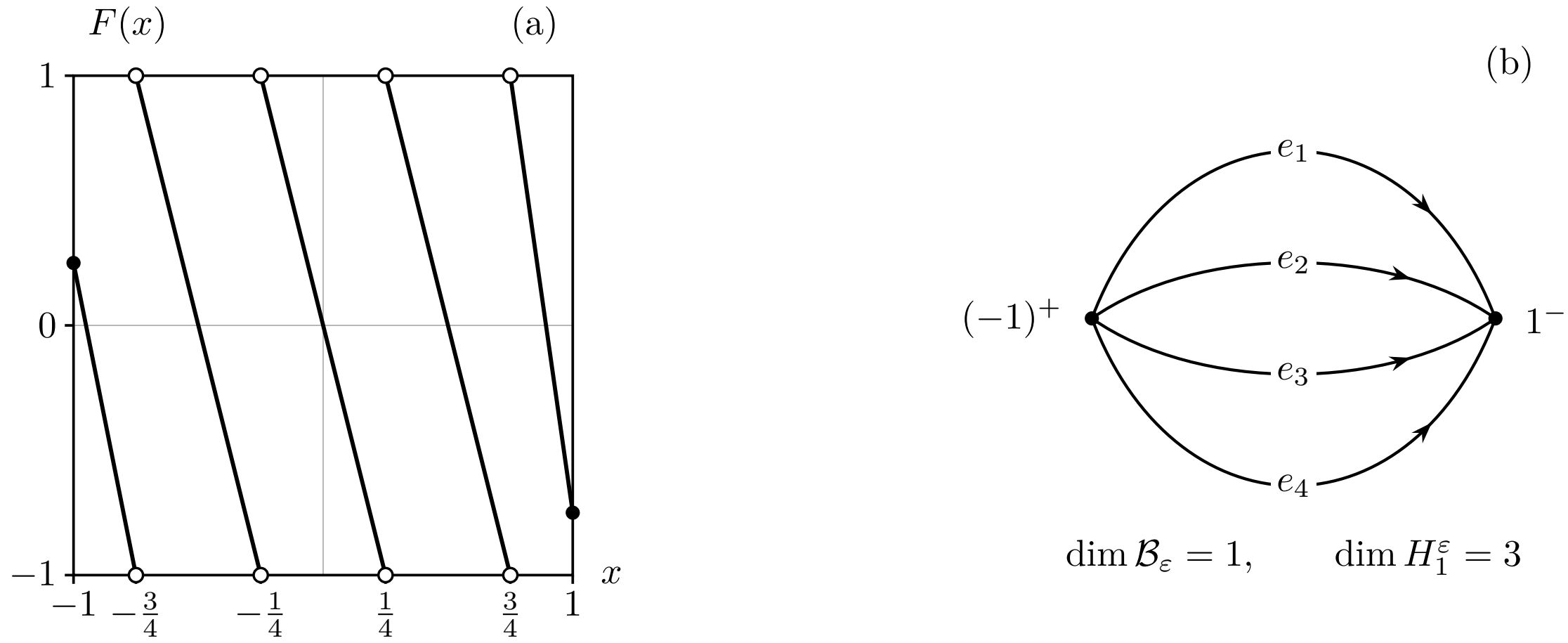

On $[-1,1]$ take cuts $-3/4,-1/4,1/4,3/4$ and parameters
$\alpha,\beta\in(-1,1)$. Define
\begin{equation}\label{eq:family}
 F_{\alpha,\beta}(x)=
 \begin{cases}
 -1-4(\alpha+1)(x+3/4),&-1\le x<-3/4,\\
 -4x-2,&-3/4<x<-1/4,\\
 -4x,&-1/4<x<1/4,\\
 -4x+2,&1/4<x<3/4,\\
 1-4(1-\beta)(x-3/4),&3/4<x\le1.
 \end{cases}
\end{equation}
All jumps have limits $-1,1$. The map is odd exactly when
$\beta=-\alpha$, but its jump rank is one throughout the parameter
square.
Figure~\ref{fig:rank-one} shows an asymmetric member and its jump graph.

For $(\alpha,\beta)=(1/4,-3/4)$, the endpoint germ itineraries are
\[
 (0,2,0,2,\ldots)\qquad\text{and}\qquad(4,1,4,1,\ldots).
\]
Hence
\[
 K_F(t)=\frac{4+t}{1-t^2},\qquad
 D_F(t)=\frac{1-4t-2t^2}{1-t^2},\qquad
 \mathscr L_F(t)=\frac{1+t}{1-4t-2t^2}.
\]
The Markov matrix on the five original branches is
\[
 M=\begin{pmatrix}
 1&1&1&0&0\\1&1&1&1&1\\1&1&1&1&1\\
 1&1&1&1&1\\0&1&1&1&1
 \end{pmatrix},
\qquad
 \det(\zeta\Id_5-M)=\zeta^2(\zeta-1)(\zeta^2-4\zeta-2).
\]
Therefore
\[
 h_{\mathrm{top}}(F)=\log(2+\sqrt6),\qquad
 \ell_0=1,\quad\ell_1=5,\quad
 \ell_n=4\ell_{n-1}+2\ell_{n-2}\quad(n\ge2).
\]

\subsection{Two different jump images and a determinant of order two}
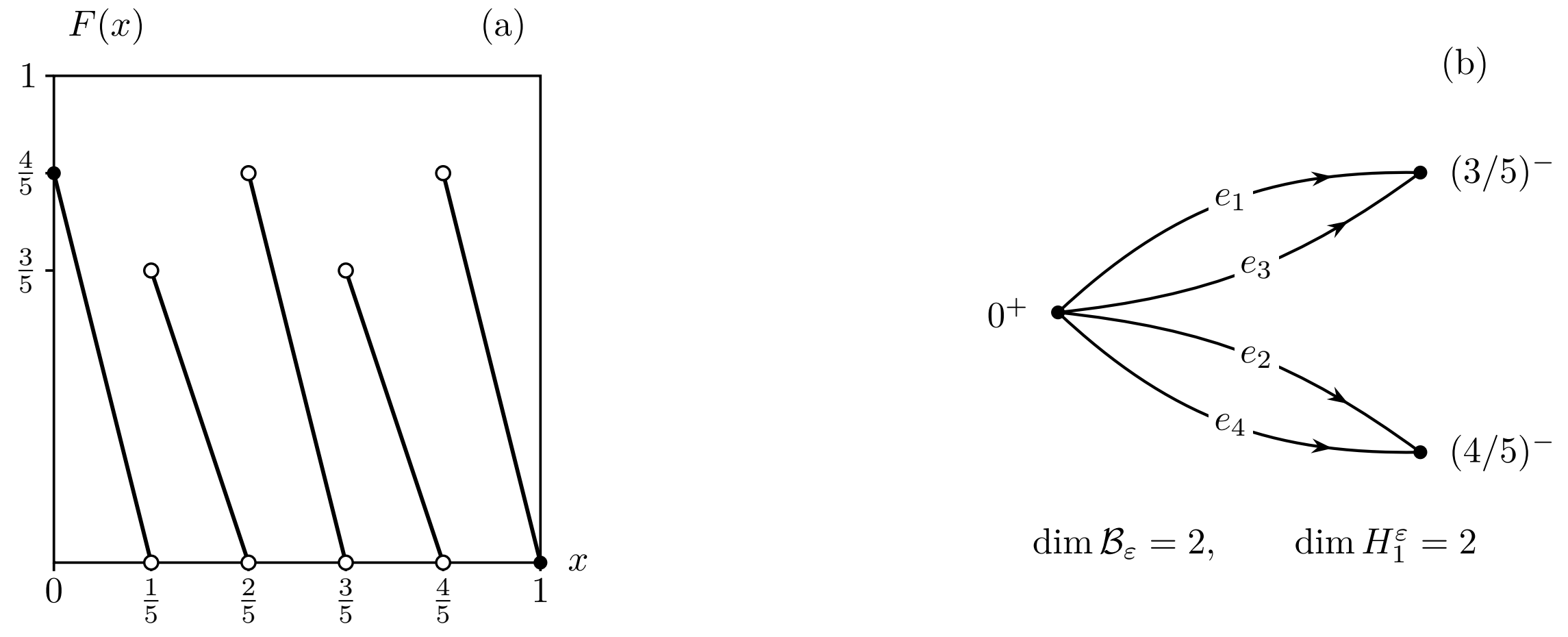
\begin{figure}[htbp]
\centering
\begin{minipage}[c]{.49\textwidth}
\centering
\begin{tikzpicture}[x=4.7cm,y=4.7cm,font=\small]
 \draw[map frame] (0,0) rectangle (1,1);
 \foreach \x/\lab in {0/{0},.2/{\frac15},.4/{\frac25},
                      .6/{\frac35},.8/{\frac45},1/{1}}{
  \draw[map frame] (\x,0)--(\x,-.0175);
  \node[below,inner sep=2pt] at (\x,-.0175) {$\lab$};
 }
 \foreach \y/\lab in {.6/{\frac35},.8/{\frac45},1/{1}}{
  \draw[map frame] (0,\y)--(-.0175,\y);
  \node[left,inner sep=2pt] at (-.0175,\y) {$\lab$};
 }
 \node[anchor=south west] at (0,1.035) {$F(x)$};
 \node[anchor=west] at (1.0275,0) {$x$};
 \node[anchor=south east] at (1,1.035) {(a)};
 \draw[map branch] (0,.8)--(.2,0);
 \draw[map branch] (.2,.6)--(.4,0);
 \draw[map branch] (.4,.8)--(.6,0);
 \draw[map branch] (.6,.6)--(.8,0);
 \draw[map branch] (.8,.8)--(1,0);
 \foreach \x in {.2,.4,.6,.8} \node[open end] at (\x,0) {};
 \foreach \x in {.2,.6} \node[open end] at (\x,.6) {};
 \foreach \x in {.4,.8} \node[open end] at (\x,.8) {};
 \node[closed end] at (0,.8) {};
 \node[closed end] at (1,0) {};
\end{tikzpicture}
\end{minipage}\hfill
\begin{minipage}[c]{.49\textwidth}
\centering
\begin{tikzpicture}[font=\small]
 \path[use as bounding box] (-.6,-2.8) rectangle (4.4,2.8);
 \node[anchor=north east] at (4.3,2.7) {(b)};
 \coordinate (z0) at (0,0);
 \coordinate (z1) at (3.5,1.35);
 \coordinate (z2) at (3.5,-1.35);
 \draw[jump edge] (z0) to[bend left=22]
  node[edge label,pos=.52] {$e_1$} (z1);
 \draw[jump edge] (z0) to[bend right=15]
  node[edge label,pos=.52] {$e_3$} (z1);
 \draw[jump edge] (z0) to[bend left=15]
  node[edge label,pos=.52] {$e_2$} (z2);
 \draw[jump edge] (z0) to[bend right=22]
  node[edge label,pos=.52] {$e_4$} (z2);
 \node[closed end] at (z0) {};
 \node[closed end] at (z1) {};
 \node[closed end] at (z2) {};
 \node[left=4pt] at (z0) {$0^+$};
 \node[right=4pt] at (z1) {$(3/5)^-$};
 \node[right=4pt] at (z2) {$(4/5)^-$};
 \node at (1.9,-2.25) {$\dim\calB_\eps=2,\qquad \dim H_1^\eps=2$};
\end{tikzpicture}
\end{minipage}
\caption{(a) The map with jump heights $u=3/5$ and $v=4/5$.
(b) The two pairs of parallel jump edges give two independent
boundaries, hence a determinant of order two.}
\label{fig:rank-two}
\end{figure}

On $[0,1]$ take cuts $c_i=i/5$, $1\le i\le4$.
All left limits are zero. The right limits are $u$ at $c_1,c_3$ and
$v$ at $c_2,c_4$, with $0<u<v<1$. Take the interior branches affine
between these prescribed limits, and the exterior branches affine
with $F(0)=v$ and $F(1)=0$. This defines a two-parameter family of
strictly decreasing branches. Its jump graph has vertices $0^+,u^-,v^-$;
its underlying simple graph is a tree, with each of the two edges
doubled. Its jump rank is two.
Figure~\ref{fig:rank-two} shows the map and graph at $u=3/5$, $v=4/5$.

At the interior parameter value $u=3/5$, $v=4/5$, put
$z_0=0^+$, $z_1=(3/5)^-$, $z_2=(4/5)^-$. Their germ transitions are
$z_0\mapsto z_2\mapsto z_0$ and $z_1\mapsto z_0$.
Let $U:\R^2\to\calB_\eps$ have columns $[z_1]-[z_0],[z_2]-[z_0]$, and write $\partial=UL$. Then
\[
 L=\begin{pmatrix}1&0&1&0\\0&1&0&1\end{pmatrix}.
\]
The columns of $K_tU$ are
\begin{align*}
 \theta(z_1,t)-\theta(z_0,t)
   &=(1,1,0,0)^{\mathsf T}
       +\frac{t}{1-t}(1,1,1,0)^{\mathsf T},\\
 \theta(z_2,t)-\theta(z_0,t)
   &=\frac{1}{1-t}(1,1,1,0)^{\mathsf T}.
\end{align*}
In the basis specified by $U$, the reduced operator has matrix
\[
 R_\eps(t)=\Id_2-\frac{t}{1-t}
 \begin{pmatrix}1+t&2\\1&1\end{pmatrix},
 \qquad D_F(t)=\frac{1-3t-2t^2}{1-t}.
\]
Here the Markov matrix is
\[
 M=\begin{pmatrix}
 1&1&1&1&0\\1&1&1&0&0\\1&1&1&1&0\\
 1&1&1&0&0\\1&1&1&1&0
 \end{pmatrix},
 \qquad \det(\zeta\Id_5-M)=\zeta^3(\zeta^2-3\zeta-2).
\]
Thus
\[
 h_{\mathrm{top}}(F)=\log\frac{3+\sqrt{17}}2,\qquad
 \mathscr L_F(t)=\frac{(1+t)^2}{1-3t-2t^2}.
\]
The determinant of order two records the two independent jump relations.

\subsection{Nonlinear branches with the same entropy}

Choose one affine branch $F_j$ of either explicit Markov map
above: $F_{1/4,-3/4}$ or the map with $u=3/5$, $v=4/5$. For a real parameter $\tau$,
replace it by
\begin{equation}\label{eq:polynomial-perturbation}
 (F_\tau)_j(x)=F_j(x)+\tau(x-c_j)^2(c_{j+1}-x)^2,
 \qquad x\in[c_j,c_{j+1}],
\end{equation}
and keep the other branches unchanged.
For sufficiently small $|\tau|$, the derivative remains negative
and the endpoint images are unchanged. The new polynomial branch
therefore has the same transition matrix and entropy.
Starting from the odd member $F_{1/4,-1/4}$ of the first family,
this perturbation on an outer branch breaks oddness while
preserving all one-sided cut itineraries.

\section{Applications to increasing branches}\label{sec:increasing}

Assume first that every branch is increasing. The jump graph is bipartite, so Theorem~\ref{thm:rank} gives $r_\eps=|\calG|-c(\Gamma_F)$. With common endpoint limits
\[
 F(c_i^-)=b,\qquad F(c_i^+)=a\qquad(1\le i\le m),
\]
one has $\partial_\eps e_i=[a^+]-[b^-]$ and $r_\eps=1$. Write
\[
 K_F(t)=\one^{\mathsf T}\xi(t)
 =\sum_{n\ge0}t^n
 \bigl(\iota(\ehat^n b^-)-\iota(\ehat^n a^+)\bigr).
\]
The rank-one calculation in Corollary~\ref{cor:rank-one} gives
\begin{equation}\label{eq:increasing-scalar}
 D_F(t)=1-tK_F(t),\qquad
 \mathscr L_F(t)=\frac1{(1-t)D_F(t)}.
\end{equation}
The entropy bounds in that corollary also hold, since every interior branch maps onto $(a,b)$. For odd maps with symmetric cuts, the summand in $K_F$ is $2\iota(\ehat^n b^-)-m$.

\subsection{Symmetric Lorenz maps and finite critical returns}

A \emph{symmetric increasing Lorenz map} is an odd map $F:[-1,1]\setminus\{0\}\to[-1,1]$, continuous and strictly increasing on each branch, with $F(0^-)=1$ and $F(0^+)=-1$. The endpoint values are their inward limits. Set $G=-F$ and let $\mathfrak r(x)=-x$. Oddness gives
\begin{equation}\label{eq:lorenz-iterate}
 G^n=\mathfrak r^{\,n}F^n
\end{equation}
wherever the orbit avoids the cut. This relation between increasing and decreasing symmetric Lorenz maps is studied in \cite{AnusicBruinCinc2020}. The two maps have the same cut preimages, hence the same cylinders and entropy.

Suppose that the germ $0^-$ first returns to the cut at time $p\ge2$. Write $\pi(x^\sigma)=x$ and put
\[
 x_j=\pi(\ehat^j1^-),\quad 0\le j\le p-1,
 \qquad x_0=1,\quad x_{p-1}=0.
\]
Define the polynomials
\begin{equation}\label{eq:lorenz-polynomials}
 P_p(t)=1-\sum_{j=0}^{p-2}\sig(x_j)t^{j+1},
 \qquad C_p(t)=1+t+\cdots+t^{p-1}.
\end{equation}
The ordered set $\{0,\pm x_0,\ldots,\pm x_{p-2}\}$ is a symmetric Markov partition. Let $M_F,M_G$ be its transition matrices and let $M_F^\pm$ be the restrictions of $M_F$ to the eigenspaces of the interval reversal $J_M$ from Section~\ref{sec:symmetry}.

\begin{theorem}[Critical-return factors]\label{thm:lorenz-factors}
The partition has $2p-2$ intervals, and
\begin{align}
 D_F(t)&=\frac{P_p(t)}{1-t^p},&
 \mathscr L_F(t)&=\frac{C_p(t)}{P_p(t)},\label{eq:lorenz-counts}\\
 \det(\Id-tM_F^+)&=P_p(t),&
 \det(\Id-tM_F^-)&=C_p(t).\label{eq:lorenz-sectors}
\end{align}
Consequently,
\[
 \det(\Id-tM_F)=P_p(t)C_p(t),\qquad
 \det(\Id-tM_G)=P_p(t)C_p(-t).
\]
\end{theorem}
\begin{proof}
No two of $x_0,\ldots,x_{p-2}$ agree or are negatives of one another: an agreement would make the orbit periodic before reaching zero, and an opposite pair would give an earlier arrival at zero by oddness. The partition therefore has $2p-1$ vertices. Side preservation gives two disjoint critical germ cycles of length $p$. The endpoint address signs are $\sig(x_0),\ldots,\sig(x_{p-2}),-1$, where the last sign records $0^-$. Hence
\[
 K_F(t)=\frac{\sum_{j=0}^{p-2}\sig(x_j)t^j-t^{p-1}}{1-t^p},
\]
and \eqref{eq:increasing-scalar} gives \eqref{eq:lorenz-counts}.

In Theorem~\ref{thm:Markov}, $\alpha=1$, the operator $Q$ is nilpotent because every interior vertex reaches the cut, and $P_\eps$ has exactly the two critical germ cycles. Thus
\[
 \det(\Id-tP_\eps)=(1-t^p)^2,
 \qquad \det(\Id-tM_F)=P_p(t)C_p(t).
\]
To identify the factors, orient each partition interval from left to right. Write $[x]$ for the basis vector of a partition vertex and $[(u,v)]$ for that of the interval between consecutive vertices $u<v$. Define the incidence boundary by
\[
 \partial_I[(u,v)]=[v]-[u].
\]
It identifies interval chains with vertex chains of total coefficient zero. For vertex reflection $R_0[x]=[-x]$, one has $\partial_IJ_M=-R_0\partial_I$. Thus antisymmetric interval chains correspond to reflection-invariant vertex chains of total coefficient zero.

The full reflection-invariant vertex space has basis
\[
 v_j=[x_j]+[-x_j]\quad(0\le j\le p-2),\qquad v_{p-1}=2[0].
\]
For $0\le j\le p-2$, let $w_j$ be the sum of the partition intervals between $0$ and $|x_j|$, minus its reflection under $J_M$. Then $\partial_Iw_j=v_j-v_{p-1}$, and the lateral images at zero give
\[
 \partial_I M_F^{\mathsf T}w_j
 =[F(|x_j|)]+[F(-|x_j|)]-[1]-[-1]
 =v_{j+1}-v_0.
\]
The vectors $v_j-v_{p-1}$ form a basis of the subspace of reflection-invariant vertex chains with total coefficient zero. Thus the induced operator is the restriction of the cyclic map $v_j\mapsto v_{j+1}$, $v_{p-1}\mapsto v_0$. Since transposition preserves sector determinants, $\det(\Id-tM_F^-)=(1-t^p)/(1-t)=C_p(t)$. The full determinant then gives the symmetric factor in \eqref{eq:lorenz-sectors}. Finally $M_G=M_FJ_M$ and $M_FJ_M=J_MM_F$, which changes the sign only on the antisymmetric sector.
\end{proof}

The antisymmetric eigenvalues are the $p$th roots of unity other than one. In particular, every eigenvalue outside the unit circle belongs to the symmetric sector.

\subsection{Periodic-orbit counts under sign reversal}

A periodic point is \emph{regular} when its orbit avoids the cut. Let $a_n(F)$ count regular fixed points of $F^n$, and let $o_n(F)$ count regular cycles of least period $n$. When these counts are finite, define the formal dynamical zeta function
\[
 \zeta_F(t)=\exp\!\left(\sum_{n\ge1}\frac{a_n(F)}n t^n\right)
           =\prod_{n\ge1}(1-t^n)^{-o_n(F)}.
\]
The cylinders of $F$ are \emph{eventually expanding} if there are $c>0$ and $\lambda>1$ such that, for every $n\ge1$, $|F^n(x)-F^n(y)|\ge c\lambda^n|x-y|$ whenever $x,y$ lie in the same $n$-cylinder. Uniform expansion means that $|F(x)-F(y)|\ge\lambda|x-y|$ on each branch for some $\lambda>1$; it implies this condition.

\begin{corollary}[Boundary correction]\label{cor:lorenz-zeta-finite}
Under the hypotheses of Theorem~\ref{thm:lorenz-factors}, if the cylinders are eventually expanding, then
\begin{equation}\label{eq:lorenz-zeta-finite}
 \zeta_F(t)=\frac{(1-t)(1-t^p)}{P_p(t)},\qquad
 \zeta_G(t)=\frac{(1+t)(1-t^p)}{P_p(t)}.
\end{equation}
\end{corollary}
\begin{proof}
Expansion makes each closed Markov path determine a unique point. Paths at partition boundaries correspond to periodic germs and must be removed. The two $F$-germ cycles have length $p$. Equation~\eqref{eq:lorenz-iterate}, applied to germs, shows that $G$ has two cycles of length $p$ when $p$ is even, and one cycle of length $2p$ when $p$ is odd. Write $\one_A=1$ when condition $A$ holds and $\one_A=0$ otherwise. Then
\[
 a_n(F)=\operatorname{tr}(M_F^n)-2p\one_{p\mid n},\qquad
 a_n(G)=\operatorname{tr}(M_G^n)-2p\one_{p\mid n,\,2\mid n}.
\]
The antisymmetric trace is $p\one_{p\mid n}-1$ for $F$, and $(-1)^n(p\one_{p\mid n}-1)$ for $G$. Subtraction gives $a_n(G)-a_n(F)=1-(-1)^n$. Using the determinant formula for $F$ and then summing this difference proves \eqref{eq:lorenz-zeta-finite}.
\end{proof}

The period-forcing results of \cite{AnusicBruinCinc2020} describe which least periods occur for symmetric Lorenz maps. Under uniform expansion, the increasing--decreasing correspondence also gives exact relations between their regular periodic-orbit counts.

Let $F$ be any symmetric increasing Lorenz map and set $G=-F$.

\begin{theorem}[Regular cycles under sign reversal]\label{thm:lorenz-zeta-general}
If $F$ is uniformly expanding and $F(1)<1$, then
\begin{equation}\label{eq:lorenz-zeta-general}
 a_n(G)-a_n(F)=1-(-1)^n,\qquad
 \frac{\zeta_G(t)}{\zeta_F(t)}=\frac{1+t}{1-t}.
\end{equation}
Consequently $o_1(G)=o_1(F)+2$, $o_2(G)=o_2(F)-1$, and $o_n(G)=o_n(F)$ for every $n\ge3$.
\end{theorem}
\begin{proof}
Let $\lambda>1$ be an expansion constant for $F$. Define $Q_0:[0,1]\to[0,1]$ by $Q_0(0)=1$ and $Q_0(u)=|F(u)|$ for $0<u\le1$. Expansion gives $0<\lambda-1\le F(1)<1$, so $Q_0$ is continuous, with one turning point $c$ and $Q_0(c)=0$. Fix $n\ge1$ and set $H(u)=Q_0^n(u)-u$. The function $H$ is positive near zero and negative near one. If an endpoint is fixed, these signs follow from expansion by at least $\lambda^n$ on its adjacent branch; otherwise they follow by continuity.

At a fixed point whose orbit avoids $c$, an increasing branch of $Q_0^n$ crosses the diagonal from below to above, and a decreasing branch crosses from above to below. A fixed point whose orbit meets $c$ belongs to the critical periodic orbit. At each of its interior points, $Q_0^n$ has a local extremum; expansion makes $H$ have the same sign on both sides. Such a point contributes no crossing. There are finitely many monotone branches and at most one fixed point per branch. Let $N_+$ count regular crossings on increasing branches and $N_-$ those on decreasing branches. The endpoint signs give $N_--N_+=1$.

Oddness gives $Q_0^n(u)=|F^n(u)|$ on regular orbits. Since each branch of $F^n$ is increasing, a positive regular fixed point of $Q_0^n$ lifts to two fixed points of $F^n$ exactly when its branch is increasing. For odd $n$, \eqref{eq:lorenz-iterate} makes the decreasing branches correspond instead to the fixed points of $G^n$. Thus $a_n(G)-a_n(F)=2$ for odd $n$. For even $n$, $G^n=F^n$. Summing gives the zeta quotient, and its factorisation $(1-t^2)/(1-t)^2$ gives the least-period counts.
\end{proof}

For example, define $L_\beta(x)=\beta x-\sig(x)$ for $x\in[-1,1]\setminus\{0\}$ and $1<\beta\le2$. At $\beta=(1+\sqrt5)/2$, the first critical return is $p=3$ and $P_p(t)=1-t-t^2$. Thus
\[
 \begin{aligned}
 \mathscr L_{L_\beta}(t)&=\frac{1+t+t^2}{1-t-t^2},\\
 \det(\Id-tM_{L_\beta})&=(1-t-t^2)(1+t+t^2),\\
 \det(\Id-tM_{-L_\beta})&=(1-t-t^2)(1-t+t^2).
 \end{aligned}
\]
Theorem~\ref{thm:lorenz-zeta-general} applies to every $1<\beta<2$. At $\beta=2$, including the endpoints among regular periodic points gives $a_n(L_2)=a_n(-L_2)=2^n$, and hence equal zeta functions.

{\color{black}
For $\beta=3/2$, a finite critical return would make $\beta$ a root of a monic polynomial with integer coefficients, contradicting the rational root theorem. Thus $L_{3/2}$ has no finite critical return, while Theorem~\ref{thm:lorenz-zeta-general} gives $\zeta_{-L_{3/2}}(t)/\zeta_{L_{3/2}}(t)=(1+t)/(1-t)$ and equal numbers of cycles of every least period at least three for $L_{3/2}$ and $-L_{3/2}$.
\par}

For $1<\beta<2$, the map $L_\beta$ has \emph{lateral matching} if the base points of $\widehat L_\beta^n0^-$ and $\widehat L_\beta^n0^+$ agree for some $n\ge1$. Let $\mathcal T\subset(1,2)$ be the set of slopes with this property.

\begin{remark}[Matching and parameter topology]\label{rem:lorenz-matching}
The set $\mathcal T$ is countable and dense in $(1,2)$. For every continuous nonconstant map $b:[0,1]\to(1,2]$, the set $b^{-1}((1,2]\setminus\mathcal T)$ is not open in $[0,1]$.
\end{remark}
\begin{proof}
Symmetry makes lateral matching equivalent to a critical return. A return after $q$ endpoint iterates gives a nonzero equation $\beta^q-\sum_{j=0}^{q-1}\sigma_j\beta^{q-1-j}=0$ with $\sigma_j\in\{-1,1\}$, proving countability; compare \cite{BruinCarminatiKalle2017}. For density, $u\mapsto|\beta u-1|$ on $[0,1]$ is conjugate to $T_\beta(v)=\min\{\beta v,2-\beta v\}$ on its invariant interval $[2-\beta,1]$ by $u=(1-v)/(\beta-1)$. Critical periodicity is dense in this tent family \cite[Corollary~4.7]{LindseyWu2023}.

If $b^{-1}((1,2]\setminus\mathcal T)$ were open, the compact set $b^{-1}(\mathcal T)$ would have image $\mathcal T\cap b([0,1])$, which is dense and not closed in a nondegenerate interval. This contradiction proves the second assertion.
\end{proof}

For the reduction proposed in \cite[Conjecture~4.3]{CholewaIgra2025}, suppose that a path of symmetric factors $L_{b(s)}$, $s\in[0,1]$, is realized with $b:[0,1]\to(1,2]$ continuous and nonconstant. By Remark~\ref{rem:lorenz-matching}, the parameters with nonmatching factors do not form an open subset of $[0,1]$. The relation of these factors to Lorenz flows is studied in \cite{CholewaIgra2026}.

\subsection{Surjective shadowing and itinerary cancellation}
\label{sec:shadowing}

In this section the map is defined on the whole closed interval, including its values at the cut and endpoints. For such a self-map $f$, a sequence $(x_j)_{j\geq0}$ is
a \emph{$\delta$-pseudo-orbit} if $|f(x_j)-x_{j+1}|<\delta$ for every $j\geq0$.
The map has \emph{shadowing} if for every $\epsilon>0$ there exists
$\delta>0$ such that each $\delta$-pseudo-orbit is
$\epsilon$-shadowed by an exact orbit $(f^j(y))_{j\geq0}$:
$|f^j(y)-x_j|<\epsilon$ for every $j\geq0$.

Shadowing for tent maps and continuous piecewise linear maps was studied by Coven, Kan and Yorke \cite{CovenKanYorke1988} and Chen \cite{Chen1991}; Portela \cite{Portela2014} considered discontinuous piecewise linear maps. Bura, Good and Samuel \cite[Section~5]{BuraGoodSamuel2026} asked whether a surjective discontinuous interval map can possess this property. We answer their question affirmatively by an explicit two-parameter family.

The construction combines an expanding branch with a forward-invariant contracting interval. This gives a surjective zero-entropy example alongside the injective piecewise contractions studied in \cite{CalderonVillar2024Contracting}.

Fix parameters $1<\mu<k<2$ and consider the piecewise affine map
\begin{equation}\label{eq:shadowmap}
 f_{k,\mu}(x)=
 \begin{cases}
  kx,&0\leq x<1,\\
  \mu+(2-\mu)(x-1),&1\leq x\leq2.
 \end{cases}
\end{equation}

Its entropy and kneading determinant are computed from the two open branches, as in Section~\ref{sec:setting}. For shadowing, the assigned value $f_{k,\mu}(1)=\mu$ is part of the map.

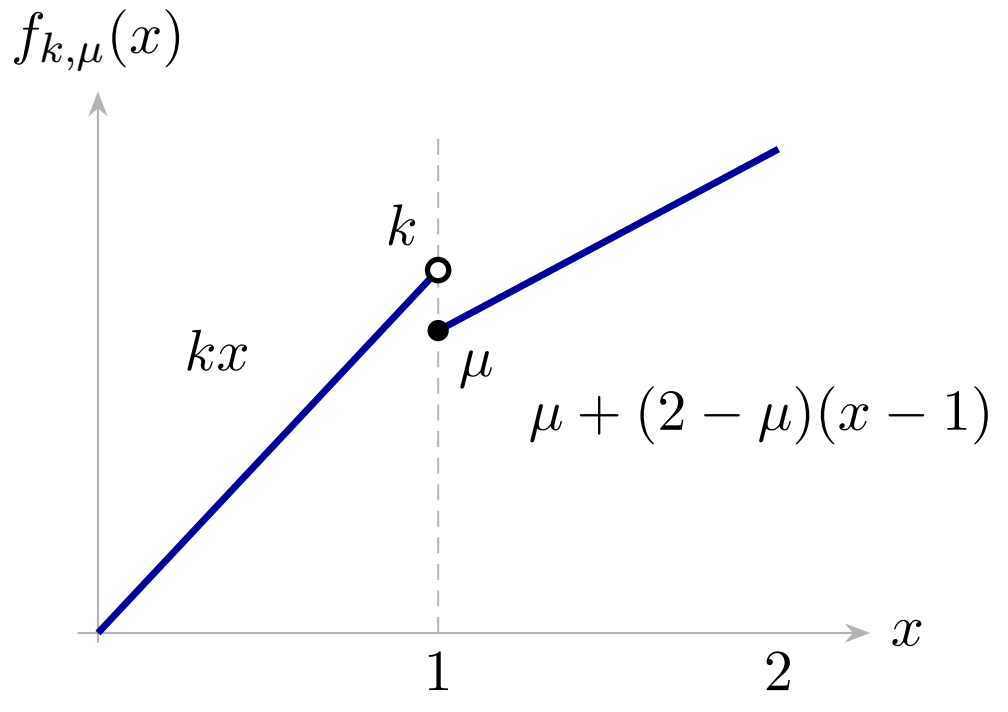
\begin{figure}[htbp]
\centering
\begin{tikzpicture}[x=2.25cm,y=1.60cm,>=Stealth]
\draw[->,gray!60] (-.06,0)--(2.27,0) node[right,black] {$x$};
\draw[->,gray!60] (0,-.04)--(0,2.24) node[above,black] {$f_{k,\mu}(x)$};
\draw[densely dashed,gray!55] (1,0)--(1,2.05);
\draw[blue!60!black,line width=1.3pt] (0,0)--(1,1.5);
\draw[blue!60!black,line width=1.3pt] (1,1.25)--(2,2);
\draw[fill=white,line width=1pt] (1,1.5) circle (2pt);
\fill (1,1.25) circle (2pt);
\node[above left] at (1,1.52) {$k$};
\node[below right] at (1,1.25) {$\mu$};
\node[below] at (1,0) {$1$};
\node[below] at (2,0) {$2$};
\node at (0.35,1.17) {$kx$};
\node at (1.95,.9) {$\mu+(2-\mu)(x-1)$};
\end{tikzpicture}
\caption{The surjective example with $(k,\mu)=(3/2,5/4)$. The left branch is expanding and the right branch is contracting. The filled point specifies the value at the cut.}
\label{fig:shadowing}
\end{figure}

\begin{theorem}[Surjective shadowing with a genuine jump]
\label{thm:surjective-shadowing}
Every map $f_{k,\mu}$ in \eqref{eq:shadowmap} is surjective and discontinuous,
is strictly increasing on its two branches, and has shadowing.
Its topological entropy is zero and its interior kneading determinant is identically one.
\end{theorem}
\begin{proof}
The two branch images are $[0,k)$ and $[\mu,2]$. Their union is $[0,2]$
because $\mu<k$, and the jump at $x=1$ has size $k-\mu>0$.

Put $\lambda=2-\mu\in(0,1)$. Given $\epsilon>0$, choose $\delta>0$
such that
\[
\delta<\min\{\mu-1,k-1\},\qquad
\delta\left(3+\frac{1}{k-1}+\frac{1}{\mu-1}\right)<\epsilon.
\]
Let $(x_j)_{j\geq0}$ be a $\delta$-pseudo-orbit.
Since $f_{k,\mu}([1,2])=[\mu,2]$ and $\mu-\delta>1$, once a term enters
$[1,2]$ all later terms remain there.

If no term enters $[1,2]$, then $x_{j+1}>kx_j-\delta$ and $x_j<1$
for every $j$. For fixed $j$ and any $r\geq1$, iteration gives
\[
 x_j<k^{-r}x_{j+r}+\frac{\delta(1-k^{-r})}{k-1}.
\]
Letting $r\to\infty$ yields $x_j\leq\delta/(k-1)$. The fixed orbit at zero therefore
$\epsilon$-shadows $(x_j)$.

Otherwise let $N$ be the first index with $x_N\geq1$.
If $N=0$, take $y=x_0$. The contraction estimate on $[1,2]$
gives $|f_{k,\mu}^j(y)-x_j|\leq\delta/(1-\lambda)<\epsilon$ for all $j$.
If $N\geq1$, then $1\leq x_N<k+\delta$. Set
\[
 y_N=\min\{x_N,k-\delta\},\qquad
 y_j=k^{j-N}y_N\quad(0\leq j<N),
\]
so $1\leq y_N<k$, $y_j<1$ for $j<N$, and
$|x_N-y_N|<2\delta$. Thus $f_{k,\mu}(y_j)=y_{j+1}$ for $j<N$.
This choice makes the exact orbit enter the right-hand branch at the same time as the pseudo-orbit.
Writing $x_{j+1}=kx_j+e_j$ with $|e_j|<\delta$ for $j<N$ gives
\[
 |x_j-y_j|\leq 2\delta+\frac{\delta}{k-1}
 \qquad (0\leq j\leq N).
\]
Continue the exact orbit from $y_N$ in the right-hand branch.
The contraction estimate gives
\[
 |x_j-f_{k,\mu}^{\,j-N}(y_N)|
 \leq 2\delta+\frac{\delta}{1-\lambda}
 \qquad(j\geq N).
\]
Both bounds are less than $\epsilon$ by the choice of $\delta$.
This proves shadowing.

For the entropy, an admissible itinerary consists of a string of left-branch
symbols followed by right-branch symbols. Hence there are at most $n+1$
length-$n$ itineraries, and their exponential growth rate is one.
Furthermore, the two image germs of the single jump have base points
$k$ and $\mu$, both in $(1,2)$. Their signed itinerary series coincide:
$\theta(k^-,t)=\theta(\mu^+,t)=1/[2(1-t)]$.
The recurrence \eqref{eq:factor} therefore gives $D_{f_{k,\mu}}(t)=1$.
\end{proof}

\section{Concluding remarks}

A further question concerns the continuity of entropy at parameters where distinct jump-image germs coalesce and the boundary dimension changes.

\section*{Statements and declarations}

\noindent\textbf{Funding.}
This work was partially supported by Funda\c{c}\~ao para a Ci\^encia e a Tecnologia
(FCT), Portugal, through CAMGSD, under projects
UID/04459/2025 (\url{https://doi.org/10.54499/UID/04459/2025}) and
UID/PRR/04459/2025 (\url{https://doi.org/10.54499/UID/PRR/04459/2025}), financed
by national funds and by the European Union through the Recovery and Resilience
Plan (PRR)--NextGenerationEU.

\medskip\noindent\textbf{Competing interests.}
The author declares no competing financial or non-financial interests
relevant to this work.

\medskip\noindent\textbf{Data availability.}
No datasets were generated or analysed in this theoretical study.

\medskip\noindent\textbf{Author contributions.}
Henrique M. Oliveira is the sole author of this work.

\medskip\noindent\textbf{Ethics approval and consent.}
Not applicable.

\medskip\noindent\textbf{Use of generative artificial intelligence.}
Generative AI assisted with the revision of the manuscript. The author is responsible for
the final mathematical arguments and text.

\begingroup\small
\bibliographystyle{plain}
\bibliography{Homological_reduction_unified}
\endgroup
\end{document}